\documentclass[11pt]{amsart}

\usepackage{setspace}

\usepackage{longtable}
\usepackage{dcolumn}

\usepackage{amscd,amssymb}
\usepackage{amsthm,amsmath,amssymb}
\usepackage{supertabular}
\usepackage[matrix,arrow]{xy}
\newcolumntype{.}{D{.}{.}{1}}

 \makeatletter\@addtoreset{equation}{section}
\makeatother

\newcommand{\mult}{\operatorname{mult}}

\newcommand{\Aut}{\operatorname{Aut}}

\newcommand{\Bir}{\operatorname{Bir}}

\renewcommand{\emptyset}{\varnothing}

\newcommand{\ZZ}{{\mathbb Z}}

\newcommand{\QQ}{{\mathbb Q}}
\newcommand{\PP}{{\mathbb P}}

\newcommand{\FF}{{\mathbb F}}

\newcommand{\nn}{\texttt{N}}
\newcommand{\FGZero}{$\textbf{F}^0$}
\newcommand{\FGOne}{$\textbf{F}^1 $}
\newcommand{\FGTwo}{$\textbf{F}^2$}
\newcommand{\FGThree}{$\textbf{F}^3$}\newcommand{\FGFive}{$\textbf{F}^5$}
\newcommand{\EG}{$\hat{\textbf{F}}^3$}

\newtheorem{theorem}[equation]{Theorem}

\newtheorem{proposition}[equation]{Proposition}
\newtheorem{lemma}[equation]{Lemma}
\newtheorem{corollary}[equation]{Corollary}
\newtheorem{conjecture}[equation]{Conjecture}

\newtheorem{remark}[equation]{Remark}

\author{Ivan Cheltsov and Jihun Park}

\title{Weighted Fano threefold hypersurfaces}

\address{\emph{Ivan Cheltsov}\newline \textnormal{School of Mathematics, University of Edinburgh, Edinburgh EH9 3JZ, UK;  \texttt{I.Cheltsov@ed.ac.uk}}}%
\address{\emph{Jihun Park}\newline \textnormal{Department of Mathematics, POSTECH, Pohang, Kyungbuk 790-784, Korea; \texttt{wlog@postech.ac.kr}}}%

\begin{document}

\begin{abstract}
We study  birational transformations into elliptic fibrations and
birational automorphisms of quasismooth anticanonically embedded
weighted Fano 3-fold hypersurfaces with terminal singularities
classified by A.R.\,Iano-Fletcher, J.\,Johnson, J.\,Koll\'ar, and
M.\,Reid.
\end{abstract}

\maketitle

\section{Introduction.}
\label{section:introduction}

Let $S$  be a smooth cubic surface in $\PP^3$ defined over a
perfect field  $k$ with Picard rank 1. For example, the equation $
2x^{3}+3y^{3}+5z^{3}+7w^{3}=0$ defines such a cubic in
$\mathrm{Proj}\big(\mathbb{Q}[x,y,z,w]\big)$ (see \cite{Ma72} or
\cite{Se43}). The condition that the Picard rank is one simply
means that every curve on $S$ defined over $k$ is cut by some
hypersurface in $\PP^3$. The surface $S$ is proved to be
birationally rigid and hence nonrational (see \cite{Ma66}).

Let $P$ and $Q$ be distinct $k$-rational points on the surface
$S$. We then consider the projection $\phi:S\dasharrow\PP^2$ from
the point $P$. Because the map $\phi$ is a double cover
generically over $\PP^2$, it induces a birational involution
$\alpha$ of the surface $S$ that interchanges two points of a
generic fiber of the rational map $\phi$. Traditionally, the
involution $\alpha$ is called a Geiser involution.

Meanwhile, we consider the line $L\subset\PP^3$ passing through
the points  $P$ and $Q$. Then the line $L$ meets the surface $S$
at another $k$-rational point $O$. For a sufficiently general
hyperplane $H$ in $\mathbb{P}^{3}$ passing through the line $L$,
the intersection $H\cap S$ is a smooth elliptic curve $E$. Then
the reflection of the elliptic curve $E$ centered at the point $O$
induces a birational involution $\beta$ of the surface $S$ that is
called a Bertini involution.

Yu.\,Manin proved the group $\Bir(S)$ of birational automorphisms
of the surface $S$ is generated by the group $\Aut(S)$ of
biregular automorphisms and Bertini and Geiser involutions of the
surface $S$, more precisely, the sequence of groups
$$
1\to\Gamma_S\to\Bir(S)\to\Aut(S)\to 1
$$
is exact, where $\Gamma_S$ is the group generated by Bertini and
Geiser involutions. Furthermore, he also described all the
relations among these involutions (see \cite{Ma72}). These
properties mentioned so far remain true for smooth del Pezzo
surfaces of degree $1$ and $2$  with Picard rank $1$. Moreover, on
a smooth del Pezzo surface of degree $2$, the group $\Gamma_S$ is
the free product of involutions. But in the case of degree $1$,
every birational automorphism is biregular (see \cite{Ma66}).

Smooth del Pezzo surfaces of degree $1$, $2$ and $3$ are the only
smooth del Pezzo surfaces that can be anticanonically embedded
into weighted projective spaces as quasismooth hypersurfaces.
Therefore, the properties described above can be naturally
expected on anticanonically embedded quasismooth weighted Fano
3-fold hypersurfaces with terminal singularities. The first step
in this direction is done in \cite{IsMa71}, where the birational
superrigidity of smooth quartic 3-folds is proved.

Smooth quartic 3-folds are the first example of quasismooth
anticanonically embedded weighted Fano 3-fold hypersurfaces with
terminal singularities that were completely classified into 95
families by A.R.\,Iano-Fletcher, J.\,Johnson, J.\,Koll\'ar, and
M.\,Reid (see \cite{IF00} and \cite{JoKo01}) and  which were
studied quite extensively in \cite{CPR} and \cite{Ry02}.

Throughout this paper, we always let
$X\subset\mathbb{P}(1,a_{2},a_{3},a_{4},a_{5})$ be a sufficiently
general quasismooth anticanonically embedded Fano
hy\-per\-sur\-face of degree $d$ and of type $\nn$ with terminal
singularities\footnote{The weighted projective space
$\operatorname{Proj} (\mathbb{F}[x_1, x_2, \cdots, x_n])$ defined
over an arbitrary field $\mathbb{F}$
 with $\mathrm{wt}(x_i)=a_i$ be denoted by
 $\mathbb{P}_{\mathbb{F}}(a_1,  a_2,\cdots, a_n)$. The weights $a_i$
 are always assumed that $a_1\leq a_2\leq\ldots\leq a_n$. When the field
 of definition is clear, we use simply the notation $\mathbb{P}(a_1,  a_2,\cdots,
 a_n)$ instead of $\mathbb{P}_{\mathbb{F}}(a_1,  a_2,\cdots,
 a_n)$.}, where the notation $\nn$ is the entry number in
Table~\ref{table:weighted-Fanos} of Appendix.

The hypersurface $X$ is proved to be rationally con\-nected (see
\cite{Sho00}) and birationally rigid (see \cite{CPR}).
Furthermore,  it follows from \cite{CPR} that the sequence of
groups
$$
1\to \Gamma_X\to \mathrm{Bir}(X)\to \mathrm{Aut}(X)\to 1,
$$
is exact, where the group $\Gamma_X$ is a subgroup of
$\mathrm{Bir}(X)$ generated by a finite set of distinct birational
involutions $\tau_{1},\ldots,\tau_{\ell}$ explicitly described in
\cite{CPR}. All the involutions here are either an elliptic
involution or a quadratic involution. The former is a
generalization of a Bertini involution and the latter is that of
Geiser involution.

Even though the paper \cite{CPR} describes the number  of the
birational involutions $\tau_{1},\ldots,\tau_{\ell}$ and their
explicit constructions, the relations among them have been in
question. We show that the group $\Gamma_X$ has exactly one of the
following group presentations:
\[\begin{array}{lll}
\mbox{\FGZero}& =& \mbox { the trivial group};\\
\mbox{\FGOne}&=&<\tau_1 \ | \ \tau_1^2=1>;\\
\mbox{\FGTwo}&=&<\tau_1,\tau_2 \ | \ \tau_1^2=\tau^2_2=1>;\\
\mbox{\FGThree}&=&<\tau_1,\tau_2,\tau_3 \ | \
\tau_1^2=\tau^2_2=\tau^2_3=1>;\\
\mbox{\EG}&=&<\tau_1,\tau_2,\tau_3 \ | \
\tau_1^2=\tau^2_2=\tau^2_3=\tau_{1}\tau_{2}\tau_{3}\tau_{1}\tau_{2}\tau_{3}=1>;\\
\mbox{\FGFive}&=&<\tau_1,\tau_2,\tau_3, \tau_4, \tau_5 \ | \
\tau_1^2=\tau^2_2=\tau^2_3=\tau^2_4=\tau^2_5=1>,
\end{array}\]
where the generator $\tau_i$ comes from an involution of $X$ and
the group operation from the composition of maps. When the group
$\Gamma_X$ is trivial, the $3$-fold $X$ is birationally
superrigid. Also, when $X$ has a unique birational involution, the
group $\Gamma_X$ has the presentation \FGOne\ that is isomorphic
to $\ZZ/2\ZZ$. Because the number of generators of $\Gamma_X$ is
completely determined in \cite{CPR}, in order to describe the
group $\Gamma_X$, it is enough to find their relations for
$\ell\geq 2$.  We prove the following result:

\begin{theorem}
\label{theorem:birational-autmorphisms} The group $\Gamma_X$ has
the group presentation as follows:
\begin{itemize}
\item \textnormal{\FGFive}  if $\nn=7$; \item \textnormal{\EG} if
$\nn=4$, $9$, $17$, $27$; \item \textnormal{\FGThree}  if
$\nn=20$; \item \textnormal{\FGTwo}  if $\nn=5$, $6$, $12$, $13$,
$15$, $23$, $25$, $30$, $31$, $33$, $36$, $38$, $40$, $41$, $42$,
$44$, $58$, $61$, $68$, $76 $; \item \textnormal{\FGOne} if $\nn=
2$, $8$, $16$, $18$, $24$, $26$, $32$, $43$, $45$, $46$, $47$,
$48$, $54$, $56$, $60$, $65$, $69$, $74$, $79$; \item
\textnormal{\FGZero} otherwise.
\end{itemize}
\end{theorem}

This theorem with the results of \cite{CPR} can be considered as a
3-fold analogue of Yu.\,Manin's results on smooth del Pezzo
surfaces of degree $\leq 3$.

The proof of  Theorem~\ref{theorem:birational-autmorphisms} is
based on the simple observation that except the cases $\nn=7$,
$20$, $60$, the involutions $\tau_{1},\ldots,\tau_{\ell}$ are
actually elliptic and induced by a single elliptic fibration. This
shows that it is worth our while to study birational
transformations of the hypersurface $X$ into elliptic fibrations.
In particular, it is an interesting question when the 3-fold $X$
is birational to an elliptic fibration.  We prove the following
result:

\begin{theorem}
\label{theorem:elliptic-fibrations} The hypersurface
$X\subset\mathbb{P}(1,a_{2},a_{3},a_{4},a_{5})$ can be
bi\-ra\-ti\-o\-nally transformed into an elliptic fibration if and
only if $\nn\not\in \{3, 60, 75, 84, 87, 93\}$.
\end{theorem}

We remark that the hypersurface $X$ of $\nn=3$ is the only smooth
Fano 3-fold that is not birationally equivalent to an elliptic
fibration. Many examples in the 95 families of weighted Fano
3-folds have not so many ways in which we can transform them into
an elliptic fibration. Naturally, they make us expect that the
hypersurface $X$, in almost all cases, has a single birational
elliptic fibration structure (see
Conjecture~\ref{conjecture:unique-elliptic-fibration} and
Proposition~\ref{proposition:unique-elliptic-fibration}).

After the theorem above, it may be a next step to ask whether the
hypersurface $X$ can be birationally transformed to a
$K3$~fibration or not. To this question we give an affirmative
answer.

\begin{proposition}
\label{proposition:existence-of-K3-fibration} The hypersurface $X$
is birationally equivalent to a $K3$~fibration.
\end{proposition}

We should remark here that D.\,Ryder\footnote{After the early
version of this paper, he announced a paper to reinforce his
thesis. In his paper, he classified birational transformations
into $K3$ and elliptic fibrations for the cases $\nn=34, 75, 88,
90$ (see \cite{Ry05}).} has studied birational transformations of
the hypersurface $X$ into $K3$ and elliptic fibrations in his
Ph.D. thesis (see \cite{Ry02}). His thesis applied the techniques
of the papers \cite{Ch00a} and \cite{CPR} to classify all
birational transformations of $X$ to $K3$ and elliptic fibrations
in the case $\nn=5$. In addition, he constructed various kinds of
birational transformations of the hypersurface $X$ into $K3$ and
elliptic fibrations and obtained partial results on the existence
of submaximal singularities on the hypersurface $X$ in many cases.

Meanwhile, as far as we know, arithmetical properties on
quasismooth anticanonically embedded weighted Fano 3-fold
hypersurfaces with terminal singularities have never been
investigated. The papers~\cite{BoTsch99}, \cite{BoTsch00}, and
\cite{HaTsch00} give us a stimulating result that rational points
are potentially dense\footnote{The set of rational points of a
variety $V$ defined over a number field $\mathbb{F}$ is called
potentially dense if for some finite field extension $\mathbb{K}$
of the field $\mathbb{F}$ the set of $\mathbb{K}$-rational points
of the variety $V$ is Zariski dense.} on smooth Fano 3-folds
possibly except double covers of $\mathbb{P}^{3}$ ramified along
smooth sextic surfaces. In the case $\nn=1$, the potential density
of rational points on the hypersurface $X$ is proved in
\cite{HaTsch00}. The hypersurface $X$ of $\nn=2$ is birational to
a double cover of $\mathbb{P}^{3}$ ramified along a sextic surface
with $15$ nodes, which implies the potential density of rational
points (see \cite{ChPa04}). Furthermore, we prove the following:

\begin{proposition}
\label{proposition:potential-density} Suppose that $X$ is defined
over a number field. Then rational points are poten\-ti\-ally
dense on the hypersurface $X$ for $\nn=1$, $2$, $4$, $5$, $6$,
$7$, $9$, $11$, $12$, $13$, $15$, $17$, $19$, $20$, $23$, $25$,
$27$, $30$, $31$, $33$, $36$, $38$, $40$, $41$, $42$, $44$, $58$,
$61$, $68$, $76$.
\end{proposition}

It immediately follows from
Theorem~\ref{theorem:birational-autmorphisms} that the group
$\Gamma_X$ is infinite if $\ell>1$. In this case, the
constructions of the involutions $\tau_{1},\ldots,\tau_{\ell}$
easily imply that the hypersurface $X$ contains an infinitely many
rational surfaces, which implies
Proposition~\ref{proposition:potential-density} except the cases
$\nn=1, 2, 11, 19$.

Even though our main result is
Theorem~\ref{theorem:birational-autmorphisms}, for the convenience
this paper starts with the problem on existence of birational
transformations of the hypersurface $X$ into elliptic fibrations.
In Section~\ref{section:elliptic-fibrations}, we prove
Theorem~\ref{theorem:elliptic-fibrations} and classify birational
transformations of the hypersurface $X$ into elliptic fibrations
in some cases. And then
Proposition~\ref{proposition:existence-of-K3-fibration} is proved
in Section~\ref{section:fibrations-into-K3-surfaces}. We prove
Theorem~\ref{theorem:birational-autmorphisms} in
Section~\ref{section:birational-automorphisms}. Finally, we
complete the proof of
Proposition~\ref{proposition:potential-density} by proving the
potential density of rational points on $X$ in the cases $\nn=11$
and $\nn=19$.

\vspace{2mm} {\itshape Acknowledgments}. We would like to thank
F.\,Bogomolov, A.\,Borisov, A.\,Corti, M.\,Gri\-nen\-ko,
V.\,Iskov\-s\-kikh, Yu.\,Pro\-kho\-rov, V.\,Sho\-ku\-rov,
D.\,Stepanov, and M.\,Verbitsky for helpful conversations. We
thank A.\,Pukhlikov and Yu.\,Tschin\-kel for proposing us these
problems. This work was initiated when the second author visited
University of Edinburgh and we almost finished the paper while the
first author visited POSTECH in Korea. We would like to thank
University of Edinburgh and POSTECH for their hospitality. The
second author was supported by KOSEF Grant R01-2005-000-10771-0 of
Republic of Korea.

\section{Elliptic fibrations.}
\label{section:elliptic-fibrations}

In this section we prove
Theorem~\ref{theorem:elliptic-fibrations}. We start with the
simple results below that are useful for this section.

\begin{lemma}
\label{lemma:mobile-log-pairs-linear-systems-not-composed-from-pensil}
Let $Y$ be a variety and $\mathcal{M}$ be a linear system without
fixed components on the variety $Y$. If the linear system
$\mathcal{M}$ is not composed from a pencil, then there is no
Zariski closed proper subset $\Sigma\subsetneq Y$ such that
$\mathrm{Supp}(S_{1})\cap\mathrm{Supp}(S_{2})\subset\Sigma$, where
$S_{1}$ and $S_{2}$ are sufficiently general divisors of the
linear system $\mathcal{M}$.
\end{lemma}

\begin{proof}
Suppose  there is a proper Zariski closed subset $\Sigma\subset Y$
such that the set-theoretic intersection of the sufficiently
general divisors $S_{1}$ and $S_{2}$  of the linear system
$\mathcal{M}$ is contained in the set $\Sigma$. Let
$\rho:Y\dasharrow\mathbb{P}^{n}$ be the rational map induced by
the linear system $\mathcal{M}$, where $n$ is the dimension of the
linear system $\mathcal{M}$. Then there is a commutative diagram
$$
\xymatrix{
&&W\ar@{->}[ld]_{\alpha}\ar@{->}[rd]^{\beta}&&\\%
&Y\ar@{-->}[rr]_{\rho}&&\mathbb{P}^{n},&}
$$ %
where $W$ is a smooth variety, $\alpha$ is a birational morphism,
and $\beta$ is a morphism. Let $Z$ be the image of the morphism
$\beta$. Then $\mathrm{dim}(Z)\geqslant 2$ because $\mathcal{M}$
is not composed from a pencil.

Let $\Lambda$ be a Zariski closed subvariety of the variety $W$
such that the morphism
$$
\alpha\vert_{W\setminus\Lambda}:W\setminus\Lambda\longrightarrow Y\setminus\alpha(\Lambda)%
$$
is an isomorphism, and $\Delta$ be the union of the subset
$\Lambda\subset W$ and the closure of the proper transform of the
set $\Sigma\setminus\alpha(\Lambda)$ on $W$. Then $\Delta$ is a
Zariski closed proper subset of  $W$.

Let $B_{1}$ and $B_{2}$ be general hyperplane sections of the
variety $Z$, and $D_{1}$ and $D_{2}$ be the proper transforms of
the divisors $B_{1}$ and $B_{2}$  on the variety $W$ respectively.
Then $\alpha(D_{1})$ and $\alpha(D_{2})$ are general divisors of
the linear system $\mathcal{M}$. Hence, in the set-theoretic sense
we have
\[\varnothing\ne\beta^{-1}\Big(\mathrm{Supp}(B_{1})\cap\mathrm{Supp}(B_{2})\Big)
=\mathrm{Supp}(D_{1})\cap\mathrm{Supp}(D_{2})\subset\Delta\subsetneq
W\]%
because $\mathrm{dim}(Z)\geqslant 2$. However, this set-theoretic
identity is  absurd.
\end{proof}

The following result is implied by Lemma~0.3.3 in \cite{KMM} and
Lemma~\ref{lemma:mobile-log-pairs-linear-systems-not-composed-from-pensil}.

\begin{corollary}
\label{corollary:mobile-log-pairs-linear-systems-not-composed-from-pensil-and-big-and-nef-divisors}
Let $Y$ be a three-dimensional variety with canonical
singularities. Suppose that a linear system $\mathcal{M}$ on $Y$
without fixed components is not composed from a pencil. For
sufficiently general surfaces $S_1$ and $S_2$ in the linear system
$\mathcal{M}$ and  a nef and big divisor $D$, the inequality
$D\cdot S_{1}\cdot S_{2}>0$ holds.
\end{corollary}

In addition, the proof of
Lemma~\ref{lemma:mobile-log-pairs-linear-systems-not-composed-from-pensil}
implies the following result.

\begin{lemma}
\label{lemma:mobile-log-pairs-linear-systems-not-composed-from-pensil-II}
Let $Y$ be a variety.  For linear systems $\mathcal{M}$ and
$\mathcal{D}$ on $Y$ without fixed components, if the linear
system $\mathcal{M}$ is not composed from a pencil, then there is
no Zariski closed proper subset $\Sigma\subsetneq Y$ such that
$\mathrm{Supp}(S)\cap\mathrm{Supp}(D)\subset\Sigma$, where $S$ and
$D$ are general divisors of the linear system $\mathcal{M}$ and
$\mathcal{D}$, respectively.
\end{lemma}

Before we proceed, we first observe that the following hold:
\begin{itemize}
\item for $\nn=1$, a general fiber of the projection of a smooth
quartic 3-fold $X\subset\mathbb{P}^{4}$ from a line contained in $X$ is a smooth elliptic curve;%
\item for $\nn=2$, the 3-fold $X$ is birational to a double cover
of $\mathbb{P}^{3}$ ramified along a singular nodal sextic (see
\cite{ChPa04}), which is birationally equivalent to an elliptic
fibration.
\end{itemize}

\begin{lemma}
\label{lemma:existence-of-elliptic-fibration} Suppose that
$\nn\not\in \{1, 2, 3, 7, 11, 19, 60, 75, 84, 87, 93\}$. Then a
sufficiently general fiber of the natural projection
$X\dasharrow\mathbb{P}(1,a_{2},a_{3})$ is a smooth elliptic curve.
\end{lemma}

\begin{proof}
Let $C$ be a general fiber of the projection
$X\dasharrow\mathbb{P}(1,a_{2},a_{3})$. Then $C$ is not a
ra\-ti\-o\-nal curve by \cite{CPR} but $C$ is a hypersurface of
degree $d$ in
$\mathbb{P}(1,a_{4},a_{5})\cong\mathrm{Proj}(\mathbb{C}[x_1,x_4,x_5])$,
where either $\lfloor d/a_{4}\rfloor\leqslant 3$ or $\lfloor
d/a_{4}\rfloor\leqslant 4$ and $2a_{5}\leq d<2a_5+a_4$.

Let $V\subset\mathbb{P}(1,a_{4},a_{5})$ be the open subset given
by $x_1\ne 0$. Then $V\cong\mathbb{C}^{2}$ and the affine curve
$V\cap C$ is either a cubic curve when $\lfloor
d/a_{4}\rfloor\leqslant 3$  or a double cover of $\mathbb{C}$
ramified at most four points when $\lfloor d/a_{4}\rfloor\leqslant
4$ and $2a_{5}\leq d<2a_5+a_4$. Therefore, the curve $C$ is
elliptic.
\end{proof}

\begin{remark}{\normalfont
If $\nn\not\in \{2, 7, 20, 36, 60\}$, each involution $\tau_i$
generating the group $\Gamma_X$  gives the commutative diagram
$$
\xymatrix{
&X\ar@{-->}[d]_{\psi}\ar@{-->}[rr]^{\tau_{i}}&&X\ar@{-->}[d]^{\psi}\\%
&\mathbb{P}(1,a_{2},a_{3})\ar@{=}[rr]&&\mathbb{P}(1,a_{2},a_{3}),&}
$$ %
where $\psi$ is the natural projection.}\end{remark}

\begin{lemma}
\label{lemma:existence-of-elliptic-fibration-7-11-19} Suppose that
$\nn\in\{7, 11, 19\}$. Then $X$ is birational to an elliptic
fibration.
\end{lemma}

\begin{proof}
We consider only the case $\nn=19$ because in the other cases the
proof is similar.

When $\nn=19$, the hypersurface $X$ in $\mathbb{P}(1,2,3,3,4)$ can
be given by the equation
$$
x_5f_{8}(x_1,x_2,x_3,x_4,x_5)+x_3f_{9}(x_3,x_4)+x_2f_{10}(x_1,x_2,x_3,x_4,x_5)
+x_1f_{11}(x_1,x_2,x_3,x_4,x_5)=0,
$$
where $f_{i}$ is a quasi-homogeneous polynomial of degree $i$.

Let $\mathcal{H}$ be the pencil of surfaces on $X$ cut by $\lambda
x_1^{2}+\mu x_2=0$ and $\mathcal{B}$ the pencil of surfaces cut on
$X$ by $\delta x_1^{3}+\gamma x_3=0$, where
$(\delta:\gamma)\in\mathbb{P}^{1}$ and
$(\lambda:\mu)\in\mathbb{P}^{1}$. Then $\mathcal{H}$ and
$\mathcal{B}$ give a map
$$
\rho:X\dasharrow \mathbb{P}^{1}\times\mathbb{P}^{1},
$$
which is not defined in
$\mathrm{Bs}(\mathcal{H})\cup\mathrm{Bs}(\mathcal{B})$.

Let $C$ be a general fiber of $\rho$. Then $C$ is a hypersurface
in $\mathbb{P}(1,3,4)\cong\mathrm{Proj}(\mathbb{C}[x_1,x_4,x_5])$
containing the point $(0:1:0)$. Thus, the affine piece of the
curve $C$ given by $x_1\ne 0$ is a cubic curve in
$\mathbb{C}^{2}$, but $C$ is not rational (see \cite{CPR}). Hence,
the fiber $C$ is elliptic.
\end{proof}

Therefore, we have obtain
\begin{corollary}
\label{corollary:existence-of-elliptic-fibrations} If
$\nn\not\in\{3, 60, 75, 84, 87, 93\}$, then $X$ is
bi\-ra\-ti\-o\-nal to an elliptic fibration.
\end{corollary}

To  complete the proof of
Theorem~\ref{theorem:elliptic-fibrations}, we need to show that
the 3-fold $X$ is not birationally equivalent to an elliptic
fibration when $\nn\in \{3, 60, 75, 84, 87, 93\}$. However, the
paper \cite{Ch00a} shows that the 3-fold $X$ of $\nn=3$ is not
birationally equivalent to an elliptic fibration. Therefore, it is
enough to consider the cases of $\nn=60$, $75$, $84$, $87$, $93$.
Suppose that for these five cases there are a birational map
$\rho:X\dasharrow V$ and a morphism $\nu:V\to\mathbb{P}^{2}$ such
that $V$ is smooth and a general fiber of the morphism $\nu$ is a
smooth elliptic curve. We must show that these assumptions lead us
to a contradiction.

Let $\mathcal{D}=|\nu^{*}(\mathcal{O}_{\mathbb{P}^{2}}(1))|$ and
$\mathcal{M}=\rho^{-1}(\mathcal{D})$. Then $\mathcal{M}\sim
-nK_{X}$ for some natural number $n$ because the group
$\mathrm{Cl}(X)$ is generated by $-K_{X}$ (see \cite{Do82}). An
irreducible subvariety $Z\subsetneq X$ is called a center of
canonical singularities of $(X, \frac{1}{n}\mathcal{M})$ if there
is a birational morphism $f:W\to X$ and an $f$-ex\-cep\-tional
divisor $E_{1}\subset W$ such that
$$
K_{W}+\frac{1}{n}f^{-1}(\mathcal{M})\sim_{\mathbb{Q}} f^{*}(K_{X}+\frac{1}{n}\mathcal{M})+\sum_{i=1}^{m} c_{i}E_{i},%
$$
where $E_{i}$ is an $f$-exceptional divisor, $c_{1}\leqslant 0$,
and $f(E_{1})=Z$. The exceptional divisor $E_1$ is called a
maximal singularity of the log pair $(X, \frac{1}{n}\mathcal{M})$.
The set of all centers of canonical sin\-gu\-la\-ri\-ties of the
log pair $(X, \frac{1}{n}\mathcal{M})$ is denoted by
$\mathbb{CS}(X, \frac{1}{n}\mathcal{M})$.

We first show that the set $\mathbb{CS}(X,
\frac{1}{n}\mathcal{M})$ is not empty. A member of the set,
\emph{a priori}, can be a smooth point, a singular point, or a
curve on $X$. And then we show that all these cases are excluded,
which gives us a contradiction.

In what follows, we may assume that the singularities of $(X,
\frac{1}{n}\mathcal{M})$ are canonical because $X$ is birationally
rigid  by \cite{CPR}.

\begin{proposition}
\label{proposition:Nother-Fano-inequality} The singularities of
$(X, \frac{1}{n}\mathcal{M})$ are not terminal.
\end{proposition}

\begin{proof}
Suppose that the singularities of $(X, \frac{1}{n}\mathcal{M})$
are terminal. Then $(X, \epsilon\mathcal{M})$ is terminal and
$K_{X}+\epsilon\mathcal{M}$ is ample for some rational number
$\epsilon>\frac{1}{n}$. Consider the commutative diagram
$$
\xymatrix{
&&W\ar@{->}[ld]_{\alpha}\ar@{->}[rd]^{\beta}&&\\%
&X\ar@{-->}[rr]_{\rho}&&V\ar@{->}[rr]_{\nu}&&\mathbb{P}^{2},&}%
$$
where $\alpha$ and $\beta$ are birational morphisms and $W$ is
smooth. Then we have
$$
\alpha^{*}(K_{X}+\epsilon\mathcal{M})+\sum_{j=1}^{m}a_{j}F_{j}\sim_{\mathbb{Q}}K_{W}
+\epsilon\mathcal{H}\sim_{\mathbb{Q}}\beta^{*}(K_{V}+\epsilon\mathcal{D})+\sum_{i=1}^{l}b_{i}G_{i},%
$$
where $G_{i}$ is a $\beta$-exceptional divisor, $F_{j}$ is an
$\alpha$-exceptional divisor, $a_{j}$ and $b_{i}$ are rational
numbers, and $\mathcal{H}=\alpha^{-1}(\mathcal{M})$. Let $C$ be a
general fiber of $\nu\circ\beta$. Then
$$
0<C\cdot\alpha^{*}(K_{X}+\epsilon\mathcal{M})\leqslant
C\cdot(\alpha^{*}(K_{X}+
\epsilon\mathcal{M})+\sum_{j=1}^{m}a_{j}F_{j})=\beta(C)\cdot(K_{V}+\epsilon\mathcal{D})=0%
$$
because $C$ is an elliptic curve, while the divisor
$\sum_{j=1}^{m}a_{j}F_{j}$ is effective by our assumption.
\end{proof}

Consequently,  the set of centers of canonical singularities
$\mathbb{CS}(X, \frac{1}{n}\mathcal{M})$ is not empty. However, in
the sequel we will show that it is empty.

\begin{lemma}
\label{lemma:smooth-points} The set $\mathbb{CS}(X,
\frac{1}{n}\mathcal{M})$ does not contain any smooth point of $X$
\end{lemma}

\begin{proof}
See Theorem~3.1 in \cite{Co00} and Theorem~5.6.2 in \cite{CPR}.
\end{proof}

\begin{lemma}
\label{lemma:n0-curves} The set $\mathbb{CS}(X,
\frac{1}{n}\mathcal{M})$ contains no curves on $X$.
\end{lemma}

\begin{proof}
See Lemmas~3.2 and 3.5 in \cite{Ry02}
\end{proof}

Therefore, the nonempty set $\mathbb{CS}(X,
\frac{1}{n}\mathcal{M})$ can contain only singular points of $X$.
In particular, there is a point $P\in\mathrm{Sing}(X)$ such that
$P$ is a center of canonical singularities of the log pair $(X,
\frac{1}{n}\mathcal{M})$. Let $\pi:Y\to X$ be the Kawamata blow up
at the point $P$, $E$ be the exceptional divisor of $\pi$, and
$\mathcal{B}=\pi^{-1}(\mathcal{M})$. Then
$\mathcal{B}\sim_{\mathbb{Q}} -nK_{Y}$ by \cite{Ka96}.

Suppose that $-K_Y^3<0$. Let
$\overline{\mathbb{NE}}(Y)\subset\mathbb{R}^{2}$ be the cone of
effective curves of $Y$. Then the class of $-E\cdot E$ generates
an extremal ray of the cone $\overline{\mathbb{NE}}(Y)$.

\begin{lemma}
\label{lemma:cone-of-curves} There are integer numbers $b>0$ and
$c\geqslant 0$ such that $-K_{Y}\cdot (-bK_{Y}+cE)$ is numerically
equivalent to an effective irreducible  reduced curve
$\Gamma\subset Y$ and generates an extremal ray of the cone
$\overline{\mathbb{NE}}(Y)$ different from the ray generated by
$-E\cdot E$.
\end{lemma}

\begin{proof}
See Corollary~5.4.6 in \cite{CPR}.
\end{proof}

Let $S_{1}$ and $S_{2}$ be two different surfaces in
$\mathcal{B}$. Then $S_{1}\cdot S_{2}\in
\overline{\mathbb{NE}}(Y)$ but
$$
S_{1}\cdot S_{2}\equiv n^{2}K^{2}_{Y},%
$$
which implies that the class of $S_{1}\cdot S_{2}$ generates the
extremal ray of the cone $\overline{\mathbb{NE}}(Y)$ that contains
the curve $\Gamma$. However, the support of every effective cycle
$C\in \mathbb{R}^{+}\Gamma$ is contained in
$\mathrm{Supp}(S_{1}\cdot S_{2})$ because $S_{1}\cdot \Gamma<0$
and $S_{2}\cdot\Gamma<0$. Similarly, we have
$\mathrm{Supp}(S_{1}\cdot S_{2})=\Gamma$, which contradicts
Lemma~\ref{lemma:mobile-log-pairs-linear-systems-not-composed-from-pensil}
because the linear system $\mathcal{M}$ is not composed from a
pencil.

\begin{corollary}
\label{corollary:nonexistence-of-elliptic-fibrations1} The
inequality $-K_Y^3\geqslant 0$ holds.
\end{corollary}

\begin{corollary}
\label{corollary:nonexistence-of-elliptic-fibrations} When
$\nn=75$, $84$, $87$, $93$, the hypersurface $X$ is not
bi\-ra\-ti\-o\-nally equivalent to an elliptic fibration.
\end{corollary}
\begin{proof}
 The result immediately follows from
the fact that the intersection number $-K_Y^3$ is indeed negative
if $\nn=75, 84, 87, 93$ (see \cite{CPR}).
\end{proof}

From now  we consider the case $\nn=60$.  First of all, we can
conclude that the set $\mathbb{CS}(X, \frac{1}{n}\mathcal{M})$
must consist of the unique singular point
 $O$ of type $\frac{1}{9}(1,4,5)$  on $X$ because
the
 Kawamata blow ups  at the other singular points again give us negative
 $-K_Y^3$ (see \cite{CPR}). It should be pointed out that the hypersurface $X$ can be
birationally transformed into a Fano 3-fold with canonical
singularities.

Let $\pi:Y\to X$ be the Kawamata blow up at the point $O$ and
$\mathcal{B}$ be the proper transform of the linear system
$\mathcal{M}$ on the variety $Y$. Also let $P$ and $Q$ be the
singular points of the variety $Y$ contained in the exceptional
divisor $E\cong \PP(1,4,5)$ of the morphism $\pi$ that are
quotient singularities of types $\frac{1}{4}(1,1,3)$ and
$\frac{1}{5}(1,1,4)$ respectively.

\begin{lemma}
The set $\mathbb{CS}(Y, \frac{1}{n}\mathcal{B})$ contains the
point $P$.
\end{lemma}
\begin{proof}  It follows from \cite{Ka96} that
the equivalence $\mathcal{B}\sim_{\mathbb{Q}} -nK_{Y}$ holds.
Therefore, we can use the same proof of
Proposition~\ref{proposition:Nother-Fano-inequality} with nef and
big $-K_Y$ to obtain $\mathbb{CS}(Y, \frac{1}{n}\mathcal{B})\ne
\emptyset$.

We first claim that $\mathbb{CS}(Y, \frac{1}{n}\mathcal{B})$
contains at least one of the points $P$ and $Q$. Let $L$ be the
curve on $E$ corresponding to the unique curve of the linear
system $|\mathcal{O}_{\PP(1,4,5)}(1)|$. Then the curve $L$ passes
through  the points $P$ and $Q$. Since
$\mathcal{B}\sim_{\mathbb{Q}} -nK_{Y}$ we obtain
$\mathcal{B}|_E\sim_\QQ nL$. Let $Z$ be an element of the set
$\mathbb{CS}(Y, \frac{1}{n}\mathcal{B})$.

Suppose that $Z$ be a smooth point of $Y$. It then implies
$\mult_Z\mathcal{B}>n$. Let $C$ be the curve on $E$ corresponding
to a general curve in the linear system
$|\mathcal{O}_{\PP(1,4,5)}(20)|$ passing through the point $Z$.
The curve $C$ cannot be contained in the base locus of the linear
system  $\mathcal{B}$. Therefore, we obtain a contradictory
inequality
\[n=C\cdot\mathcal{B}\geq \mult_Z(C)\mult_Z(\mathcal{B})>n.\]

Suppose that $Z$ be a curve. Then $\mult_Z(\mathcal{B})\geq n$.
Let $C$ be the curve on $E$ corresponding to a general curve in
the linear system $|\mathcal{O}_{\PP(1,4,5)}(20)|$. We then have
$$
n=C\cdot\mathcal{B}\geq \mult_Z(\mathcal{B})C\cdot Z\geq n C\cdot Z,%
$$
which implies $C\cdot Z=1$ on $E$. Hence, the curve  $Z$ must be
the curve $L$.

It follows from \cite{Ka96} that if the curve $L$ belongs to the
set $\mathbb{CS}(Y, \frac{1}{n}\mathcal{B})$, then a singular
point on the curve $L$ also belongs to the set $\mathbb{CS}(Y,
\frac{1}{n}\mathcal{B})$. It proves our claim.

For now, we suppose that the set $\mathbb{CS}(Y,
\frac{1}{n}\mathcal{B})$ contains the point $Q$.

Let $\alpha:U\to Y$ be the Kawamata blow up at the point $Q$ and
$\mathcal{D}$ be the proper transform of the linear system
$\mathcal{M}$ on the variety $U$. We then see that
$\mathcal{D}\sim_{\mathbb{Q}}-nK_{U}$. The complete linear system
$|-4K_{U}|$ is the proper transform of the pencil $|-4K_{X}|$, the
base locus of which consists of a curve $Z$ such that
$\pi\circ\alpha(Z)$ is the base curve of the pencil $|-4K_{X}|$.

Let $H$ be a sufficiently general surface of the pencil
$|-4K_{U}|$. Then the equality
$$
Z^{2}=-K_{U}^{3}=-\frac{1}{30}
$$
holds on the surface $H$ but
$\mathcal{D}\vert_{H}\sim_{\mathbb{Q}} nZ$. Therefore, it follows
that
$$
\mathrm{Supp}(D)\cap\mathrm{Supp}(H)=\mathrm{Supp}(Z),
$$
where $D$ is a  general surface of the linear system
$\mathcal{D}$, which is impossible by
Lemma~\ref{lemma:mobile-log-pairs-linear-systems-not-composed-from-pensil-II}.

Consequently, the set $\mathbb{CS}(Y, \frac{1}{n}\mathcal{B})$
contains the point $P$.
\end{proof}
 The hypersurface $X$ can be given by a quasihomogeneous equation
of degree $24$
$$
x_5^{2}x_4+x_5f_{15}(x_1,x_2,x_3,x_4)+f_{24}(x_1,x_2,x_3,x_4)=0\subset\mathbb{P}(1,4,5,6,9),
$$
where $f_{i}(x_1,x_2,x_3,x_4)$ is a quasihomogeneous polynomial of
degree $i$. Let $D$ be a general surface of in the linear system
$|-5K_{X}|$ and $S$ be the unique surface of the linear system
$|-K_{X}|$. Then $D$ is cut on $X$ by the equation
$$
\lambda x_1^{5}+\delta x_1x_2+\mu x_3=0,
$$
where $(\lambda:\delta:\mu)\in\mathbb{P}^{2}$, and  $S$ is cut by
the equation $x_1=0$. Moreover, the base locus of the linear
system $|-5K_{X}|$ consists of the single irreducible curve $C$
that is cut on the hypersurface $X$ by the equations $x_1=x_3=0$.
In particular, we have $D\cdot S=C$.

In a neighborhood of the point $O$ the monomials $x_1$, $x_2$, and
$x_3$ can be considered as weighted local coordinates on $X$ such
that $\mathrm{wt}(x_1)=1$, $\mathrm{wt}(x_2)=4$ and
$\mathrm{wt}(x_3)=5$. In a neighborhood of the point $P$ the
birational morphism $\pi$ can be given by the equations
$$
x_1=\tilde{x}_{1}\tilde{x}_{2}^{\frac{1}{9}},\ x_2=\tilde{x}_{2}^{\frac{4}{9}},\
x_3=\tilde{x}_{3}\tilde{x}_{2}^{\frac{5}{9}},%
$$
where $\tilde{x}_1$, $\tilde{x}_2$ and $\tilde{x}_3$ are weighted
local coordinates on the variety $Y$ in a neighborhood of the
point $P$ such that $\mathrm{wt}(\tilde{x}_1)=1$,
$\mathrm{wt}(\tilde{x}_2)=3$ and $\mathrm{wt}(\tilde{x}_3)=1$. Let
$\tilde{D}$, $\tilde{S}$, and $\tilde{C}$ be the proper transforms
of the surface $D$, the surface $S$, and the curve $C$ on the
variety $Y$ respectively, and $E$ be the exceptional divisor of
$\pi$. Then in a neighborhood of $P$ the surface $E$ is given by
the equation $\tilde{x}_2=0$, the surface $\tilde{D}$ is given by
the equation
$$
\lambda \tilde{x}_1^{5}+\delta \tilde{x}_1+\mu \tilde{x}_3=0,
$$
and the surface $\tilde{S}$ is given by the equation
$\tilde{x}_1=0$. Hence, we see that
$$
\tilde{D}\sim_{\mathbb{Q}}\pi^{*}(-5K_{X})-\frac{5}{9}E\sim_{\mathbb{Q}}5\tilde{S}\sim_{\mathbb{Q}}-5K_{Y},
$$
the curve $\tilde{C}$ is the intersection of the surfaces
$\tilde{D}$ and $\tilde{S}$; the linear system $|-5K_{Y}|$ is the
proper transform of $|-5K_{X}|$; the base locus of $|-5K_{Y}|$
consists of the curve $\tilde{C}$.

Let $\beta:W\to Y$ be the Kawamata blow up of the point $P$. And
let $\bar{D}$, $\bar{S}$, and $\bar{C}$ be the proper transforms
on the variety $W$ of the surface $D$, the surface $S$, and the
curve $C$ respectively and $F$ be the exceptional divisor of the
morphism $\beta$. Then the surface $F$ is the weighted projective
space $\mathbb{P}(1,1,3)$ and in a neighborhood of the singular
point of the surface $F$ the birational morphism $\beta$ can be
given by the equations
$$
\tilde{x}_1=\bar{x}_1\bar{x}_2^{\frac{1}{4}},\ \tilde{x}_2=
\bar{x}_2^{\frac{3}{4}},\ \tilde{x}_3=\bar{x}_3\bar{x}_2^{\frac{1}{4}},%
$$
where $\bar{x}_1$, $\bar{x}_2$ and $\bar{x}_3$ are weighted local
coordinates on the variety $W$ in a neighborhood of the singular
point of $F$ such that $\mathrm{wt}(\bar{x}_1)=1$,
$\mathrm{wt}(\bar{x}_2)=2$ and $\mathrm{wt}(\bar{x}_3)=1$. In
particular, the exceptional divisor $F$ is given by the equation
$\bar{x}_2=0$, the surface $\bar{D}$ is given by the equation
$$
\lambda \bar{x}_1^{5}\bar{x}_2+\delta \bar{x}_1+\mu \bar{x}_3=0,
$$
and the surface $\bar{S}$ is given by the equation $\bar{x}_1=0$.
Therefore,
$$
\bar{D}\sim_{\mathbb{Q}}\beta^{*}(\tilde{D})-\frac{1}{4}F\sim_{\mathbb{Q}}
(\pi\circ\beta)^{*}(-5K_{X})-\frac{5}{9}\beta^{*}(E)-\frac{1}{4}F,\
\bar{S}\sim_{\mathbb{Q}}\beta^{*}(\bar{S})-\frac{1}{4}F\sim_{\mathbb{Q}}-K_{W},
$$
and the curve $\bar{C}$ is the intersection of the surfaces
$\bar{D}$ and $\bar{S}$. Let $\mathcal{P}$ be the proper transform
of the linear system $|-5K_{X}|$ on  $W$. Then $\bar{D}$ is a
general surface of  $\mathcal{P}$, the base locus of the linear
system $\mathcal{P}$ consists of the curve $\bar{C}$, and the
equalities
$$
\bar{D}\cdot\bar{C}=\bar{D}\cdot\bar{D}\cdot\bar{S}=\frac{1}{3}
$$
hold. Thus, the divisor $\bar{D}$ is nef and big because
$\bar{D}^{3}=2$.

Let $B_{1}$ and $B_{2}$ be general divisors of $\mathcal{D}$. Then
$$
\bar{D}\cdot B_{1}\cdot B_{2}=\Big(\beta^{*}(-5K_{Y})-\frac{1}{4}F\Big)
\cdot\Big(\beta^{*}(-nK_{Y})-\frac{n}{4}F\Big)^{2}=0,%
$$
which contradicts
Lemma~\ref{corollary:mobile-log-pairs-linear-systems-not-composed-from-pensil-and-big-and-nef-divisors}.
Hence, we have proved Theorem~\ref{theorem:elliptic-fibrations}.

One can easily check that the hypersurface $X$ can be birationally
transformed into elliptic fibrations in several distinct ways in
the case when $\nn\in\Omega=\{$$1$, $2$, $7$, $9$, $11$, $17$,
$19$, $20$, $26$, $30$, $36$, $44$, $49$, $51$, $64\}$. In other
words, in the case when $\nn\in\Omega$ there are rational maps
$\alpha:X\dasharrow\mathbb{P}^{2}$ and
$\beta\dasharrow\mathbb{P}^{2}$ such that the normalizations of
general fibers of $\alpha$ and $\beta$ are  elliptic curves but
they cannot make the diagram
$$
\xymatrix{
&X\ar@{-->}[d]_{\alpha}\ar@{-->}[rr]^{\sigma}&&X\ar@{-->}[d]^{\beta}\\%
&\mathbb{P}^{2}\ar@{-->}[rr]_{\zeta}&&\mathbb{P}^{2},&}
$$ %
commute for any birational maps $\sigma$ and $\zeta$.

\begin{conjecture}
\label{conjecture:unique-elliptic-fibration} Let
$\rho:X\dasharrow\mathbb{P}^{2}$ be a ra\-ti\-o\-nal map such that
the normalization of a general fiber of $\rho$ is an elliptic
curve. Then there is a commutative diagram
$$
\xymatrix{
&&&&X\ar@{-->}[dll]_{\psi}\ar@{-->}[dr]^{\rho}&&&\\%
&&\mathbb{P}(1,a_{2},a_{3})\ar@{-->}[rrr]_{\phi}&&&\mathbb{P}^{2},&&}
$$
if  $\nn\not\in\{3, 60, 75, 84, 87, 93\}\cup\Omega$, where $\psi$
is the natural projection and $\phi$ is a birational map.
\end{conjecture}

In the case $\nn=5$,
Conjecture~\ref{conjecture:unique-elliptic-fibration} has been
verified in \cite{Ry02}.

\begin{proposition}
\label{proposition:unique-elliptic-fibration}
Conjecture~\ref{conjecture:unique-elliptic-fibration} holds for
$\nn=14$, $22$, $28$, $34$, $37$, $39$, $52$, $53$, $57$, $59$,
$66$,  $70$, $72$, $73$, $78$, $81$, $86$, $88$, $89$, $90$, $92$,
$94$, $95$.
\end{proposition}

\begin{proof}
In the proof of Theorem~\ref{theorem:elliptic-fibrations}, we see
that there is a point $P\in\mathrm{Sing}(X)$ that belongs to
$\mathbb{CS}(X, \frac{1}{n}\mathcal{M})$. Let $\pi:Y\to X$ be the
Kawamata blow up of the point $P$, $E$ be the exceptional divisor
of $\pi$, and $\mathcal{B}$ be the proper transform on $Y$ of
$\mathcal{M}$. Then $\mathcal{B}\sim_{\mathbb{Q}} -nK_{Y}$  by
\cite{Ka96}.

There is exactly one singular point of the hypersurface $X$, say
the point $Q$, such that we have $-K_{Y}^{3}=0$ if $P=Q$, and
$-K_{Y}^{3}<0$ if $P\ne Q$. In the case when $P\ne Q$ we can
proceed as in the proof of
Theorem~\ref{theorem:elliptic-fibrations} to derive a
contradiction. Thus, we have $P=Q$.

The linear system $|-rK_{Y}|$ is free for some $r\in\mathbb{N}$
and induces a morphism
$$
\phi:Y\to \mathbb{P}(1,a_{2},a_{3})
$$
such that $\phi=\psi\circ\pi$. However, for a general surface
$S\in\mathcal{B}$ and a general fiber $C$ of the morphism $\phi$
we have $S\cdot C=0$. Hence, $\mathcal{B}$ lies in the fibers of
the elliptic fibration $\phi$, which implies the claim.
\end{proof}

Therefore,  in many cases, the hypersurface $X$ can be
birationally transformed into an elliptic fibration in a unique
way.

\section{Fibrations into $K3$ surfaces.}
\label{section:fibrations-into-K3-surfaces}

 In this section, we prove
Proposition~\ref{proposition:existence-of-K3-fibration}.

\begin{lemma}
\label{lemma:K3-fibration-18-28} Suppose that $\nn\in\{18, 22,
28\}$. Then $X$ is birational to a $K3$ fibration.
\end{lemma}

\begin{proof}
Let $\mathcal{H}$ be the pencil in $|-a_{3}K_{X}|$ of surfaces
passing through the singular point of the hypersurface $X$ of type
${\frac{1}{a_{3}}}(1,-1,1)$. Then a general surface in
$\mathcal{H}$ is a compactification of a quartic in
$\mathbb{C}^{3}$, which implies that $X$ is birational to a $K3$
fibration.
\end{proof}

Suppose that $\nn\not\in\{18, 22, 28\}$. Let
$\psi:X\dasharrow\mathbb{P}^{1}$ be the map induced by the
projection
$$
\mathbb{P}(1,a_{2},a_{3},a_{4},a_{5})\dasharrow\mathbb{P}(1,a_{2})
$$
and $S$ be a general fiber of $\psi$. Then the surface $S$ is a
hypersurface of degree $d$ in $\mathbb{P}(1,a_{3},a_{4},a_{5})$
that is not uniruled because $X$ is birationally rigid by
\cite{CPR}. Therefore, we may assume in the following that
$a_{2}\ne 1$. Let us show that $S$ is birational to a $K3$
surface.

\begin{lemma}
\label{lemma:K3-fibration-I} Suppose that $\lfloor
d/a_{3}\rfloor\leqslant 4$. Then $S$ is birational to a $K3$
surface.
\end{lemma}

\begin{proof}
The surface $S$ is a compactification of a quartic in
$\mathbb{C}^{3}$.
\end{proof}

\begin{lemma}
\label{lemma:K3-fibration-II} Suppose that $2a_{5}+a_{3}>d$ and
$\lfloor d/a_{3}\rfloor\leqslant 6$. Then the surface $S$ is
birationally equivalent to a $K3$ surface.
\end{lemma}

\begin{proof}
The surface $S$ is a compactification of a double cover of
$\mathbb{C}^{2}$ ramified along a sextic curve, which implies that
$S$ is birational to a $K3$ surface.
\end{proof}


\begin{lemma}
\label{lemma:K3-fibration-III} Suppose that  $2a_{5}+2a_{3}>d$,
$3a_5>d$, and $d\leqslant 5a_{3}$. Then the surface $S$ is
birationally equivalent to a $K3$ surface.
\end{lemma}

\begin{proof}
The surface $S$ is a compactification of a double cover of
$\mathbb{C}^{2}\setminus L$ ramified along a quintic curve, where
$L$ is a line in $\mathbb{C}^{2}$, which implies the statement.
\end{proof}

Consequently, we may consider the 3-fold $X$ only when
$$
\nn\in \{27, 33, 48, 55, 56, 58, 63, 65, 68, 72, 74, 79, 80, 83, 85, 89, 90, 91, 92, 94, 95\}.%
$$

\begin{lemma}
\label{lemma:K3-fibration-rest} Suppose that $\nn\not\in\{27, 56,
65, 68, 83\}$. Then the surface $S$ is birationally equivalent to
a $K3$ surface.
\end{lemma}

\begin{proof}
In the case $\nn=91$, the rational map $\psi$ is studied in
Example~2.5 in \cite{Ry02}, which implies that the surface $S$ is
birational to a $K3$ surface. We use the same approach for the
others. We consider only the case $\nn=72$, because the proof is
similar in other cases.

Let $X$ be a general hypersurface  in $ \mathbb{P}(1,2,3,10,15) $
of degree $30$. Let $\Gamma$ be the curve on the hypersurface $X$
given by the equation $x_1=x_2=0$ and $\mathcal{B}$ be the pencil
of surfaces on the hypersurface $X$ that are cut by the equations
$$
\lambda x_1^{2}+\mu x_2=0,%
$$
where $(\lambda:\mu)\in\mathbb{P}^{1}$. Then $S$ belongs to
$\mathcal{B}$, the curve $C$ is the base locus of the pencil
$\mathcal{B}$, and the projection $\psi$ is the rational map given
by $\mathcal{B}$. Moreover, it follows from the generality of the
hypersurface $X$ that the curve $\Gamma$ is reduced, irreducible,
and rational.

Let $P$ be a singular point of $X$ of type $\frac{1}{3}(1,2,1)$
and  $\pi:V\to X$ be the Kawamata blow up at the point $P$ with
the exceptional divisor $E\cong\mathbb{P}(1,1,2)$. Let
$\mathcal{M}$, $\hat\Gamma$, $\hat{S}$, and $\hat{Y}$  be  the
proper transforms on $V$ of the pencil $\mathcal{B}$ , the curve
$\Gamma$, the fiber $S$, and the surface $Y$ cut by the equation
$x_1=0$ on the hypersurface $X$, respectively. Then
$$
-4K_{V}^{3}=\hat{S}\cdot\hat{\Gamma}<0,
$$
where $\hat{S}\in\mathcal{M}$, the curve $\hat{\Gamma}$ is the
base locus of the pencil $\mathcal{M}$, and the equivalences
$$
\hat{S}\sim 2\hat{Y}\sim -2K_{V}\sim_{\mathbb{Q}}\pi^{*}(-2K_{X})-{\frac{2}{3}}E%
$$
hold (see Proposition~3.4.6 in \cite{CPR}). The surface $\hat{Y}$
has canonical singularities.

Let $\overline{\mathbb{NE}}(V)\subset\mathbb{R}^{2}$ be the cone
of effective curves of $V$. Then the class of $-E\cdot E$
generates one extremal ray of the cone
$\overline{\mathbb{NE}}(V)$, while the curve $\hat{\Gamma}$
generates another extremal ray of the cone
$\overline{\mathbb{NE}}(V)$ because $\hat{S}\cdot\hat{\Gamma}<0$
and $\hat{\Gamma}$ is the only base curve of the pencil
$\mathcal{M}$, which implies that the curve $\hat{\Gamma}$ is the
only curve contained in the extremal ray generated by
$\hat{\Gamma}$.

The log pair $(V, \hat{Y})$ has log terminal singularities by
Theorem~17.4 in \cite{Ko91}, which implies that the singularities
of $(V, \hat{Y})$ are canonical because $\hat{Y}\sim -K_{V}$.
Hence, for a sufficiently small rational number $\epsilon>1$ the
singularities of the log pair $(V, \epsilon\hat{Y})$ are still log
terminal but the inequality
$(K_{V}+\epsilon\hat{Y})\cdot\hat{\Gamma}<0$ holds. There is a log
flip $\alpha:V\dasharrow U$ along the curve $\hat{\Gamma}$ by
\cite{Sho93}.

Let $\mathcal{P}=\alpha(\mathcal{M})$, $\bar{Y}=\alpha(\hat{Y})$,
$\bar{S}=\alpha(\hat{S})$, and $\bar{\Gamma}$ be the flipped curve
on $U$, namely, a possibly reducible curve such that
$V\setminus\hat{\Gamma}\cong U\setminus\bar{\Gamma}$. Then the
surface  $\bar{S}$ is a member of the pencil $\mathcal{P}$, the
log pair $(U, \epsilon\bar{Y})$ has log terminal singularities,
$\bar{S}\cdot\bar{\Gamma}=2\bar{Y}\cdot\bar{\Gamma}<0$, and the
equivalences $-K_{U}\sim\bar{Y}$ and $\bar{S}\sim -2K_{U}$ hold.
Therefore, the log pair $(U, \bar{Y})$ has canonical
singularities. In particular, the singularities of the variety $U$
are canonical.

Suppose $\mathrm{Bs}(\mathcal{P})\ne\varnothing$. Then
$\mathrm{Bs}(\mathcal{P})$ consists of a possibly reducible curve
$Z$ that is numerically equivalent to $\bar{\Gamma}$. Hence, every
surface in $\mathcal{P}$ is nef. Let $H$ be a general very ample
divisor on  $V$ and $\bar{H}=\alpha(H)$. Then $\bar{H}\cdot Z<0$,
which implies $Z\subset\bar{H}$. The inequality
$$
\bar{H}\cdot\bar{S}_{1}\cdot\bar{S}_{2}<0
$$
holds for general surfaces $\bar{S}_{1}$ and $\bar{S}_{2}$ in
$\mathcal{P}$, which contradicts the numerical effectiveness of
the surface $S_{2}$ because $\bar{H}\cdot\bar{S}_{1}$ is
effective. Consequently, the pencil $\mathcal{P}$ has no base
points, and hence the surface $\bar{S}$ has canonical
singularities.

Let $\phi:U\to\mathbb{P}^{1}$ be the morphism given by the pencil
$\mathcal{P}$. Then $\bar{S}$ is a sufficiently general fiber of
$\phi$ and $2\bar{Y}$ is a fiber of $\phi$. Moreover, we have
$K_{\bar{S}}\sim 0$ by the adjunction formula because the
equivalences $-K_{U}\sim\bar{Y}$ and $\bar{Y}\vert_{\bar{S}}\sim
0$ hold. Therefore, the surface $\bar{S}$ is either an abelian
surface or a $K3$ surface.

Let $C=E\cap\hat{S}$. Then $\bar{S}$ contains $\alpha(C)$ because
$C\ne\hat{\Gamma}$ and $\alpha$ is an isomorphism in the outside
of $\hat{\Gamma}$. However, a component of $C$ must be rational
because $C$ is a hypersurface of degree $2$ in
$\mathbb{P}(1,1,2)$, which implies that $\bar{S}$ cannot be an
abelian surface\footnote{It seems to us that there are no
rationally connected 3-folds fibred into abelian or bielliptic
surfaces.}.
\end{proof}

Therefore, it is enough to check the cases $\nn\in\{27, 56, 65,
68, 83\}$ to conclude the proof of
Proposition~\ref{proposition:existence-of-K3-fibration}. We prove
that $S$ is birational to a $K3$ surface case by case.

\emph{Case $\nn=27$ or $65$.}

Because the methods for $\nn= 27$ and $65$ are the same, we only
consider the case $\nn=27$.

The surface
$S\subset\mathrm{Proj}(\mathbb{C}[x_1,x_3,x_4,x_5])\cong\mathbb{P}(1,3,5,5)$
can be given by the equation
$$
x_5^{2}f_{5}(x_1,x_3,x_4)+x_5f_{10}(x_1,x_3,x_4)+f_{15}(x_1,x_3,x_4)=0,
$$
where $f_{i}$ is a quasi-homogeneous polynomial. Introducing a
variable $y=x_5f_{5}(x_1,x_3,x_4)$ of weight $10$, we obtain the
hypersurface
$$
\tilde{S}\subset\mathbb{P}(1,3,5,10)\cong\mathrm{Proj}(\mathbb{C}[x_1,x_3,x_4,y])
$$
of degree $20$ given by the equation
$$
y^{2}+y
f_{10}(x_1,x_3,x_4)f_{5}(x_1,x_3,x_4)+f_{15}(x_1,x_3,x_4)f_{5}(x_1,x_3,x_4)=0
$$
and birational to $S$. The surface $\tilde{S}$ is a
compactification of a double cover of $\mathbb{C}^{2}$ ramified
along a sextic curve. Therefore, the surface $S$ is birational to
a $K3$ surface.

\emph{Case $\nn=56$.}

The surface $S$ is a hypersurface of degree $24$ in
$\mathbb{P}(1,3,8,11)$ given by the equation
$$
x_5^{2}x_1^{2}+x_5x_1f_{12}(x_1,x_3,x_4)+f_{24}(x_1,x_3,x_4)=0
\subset\mathrm{Proj}(\mathbb{C}[x_1,x_3,x_4,x_5]),
$$
where $f_{i}$ is a quasi-homogeneous polynomial of degree $i$.
Introducing a new variable $y=x_1x_5$ of weight $12$, we obtain
the hypersurface $\tilde{S}$ of degree $24$ in
$\mathbb{P}(1,3,8,12)$ given by the equation
$$
y^{2}+y
f_{12}(x_1,x_3,x_4)+f_{24}(x_1,x_3,x_4)=0\subset\mathrm{Proj}(\mathbb{C}[x_1,x_3,x_4,y])
$$
which is birational to $S$. We have $K_{\tilde{S}}\sim 0$, which
implies the claim.

\emph{Case $\nn=68$.}

The surface $S$ is a general quasismooth hypersurface of degree
$28$ in $\mathrm{Proj}(\mathbb{C}[x_1,x_3,x_4,x_5])$, where
$\mathrm{wt}(x_1)=1$, $\mathrm{wt}(x_3)=4$, $\mathrm{wt}(x_4)=7$,
$\mathrm{wt}(x_5)=14$. The surface $S$ has a canonical singular
point $Q$ of type $\mathbb{A}_{1}$ and two singular points $P_{1}$
and $P_{2}$ of type $\frac{1}{7}(1,4)$.

Let $\mathcal{P}$ be the pencil of curves on $S$ given by
$$
\lambda x_1^{4}+\mu x_3=0,
$$
where $(\lambda:\mu)\in\mathbb{P}^{1}$. Then the pencil
$\mathcal{P}$ gives a rational map
$\phi:S\dasharrow\mathbb{P}^{1}$ whose general fiber is an
elliptic curve. Let $\tau:Y\to S$ be the minimal resolution of
singularities, $Z$ be the proper transform on the surface $Y$ of
the irreducible curve that is cut on the surface $S$ by the
equation $x_1=0$, and $\psi=\phi\circ\tau$. Then $\psi$ is a
morphism and $Z$ lies in a fiber of $\psi$.

Consider $\tau$-exceptional curves $E$, $\hat{E}_{1}$,
$\check{E}_{1}$, $\hat{E}_{2}$, and $\check{E}_{2}$, where $
\tau(E)=Q,\ \tau(\hat{E}_{i})=P_i$, $\tau(\check{E}_{i})=P_{i},\ \check{E}_{i}^{2}=-4,%
$ and $E^{2}=\hat{E}_i^{2}=-2$. Let $L$ be the fiber of $\psi$
over the point $\psi(Z)$. Then $Z\cong\tau(Z)\cong\mathbb{P}^{1}$,
the curve $Z$ is a component of $L$ of multiplicity $4$, the fiber
$L$ contains the curve $E$, and either the surface $Y$ is a
minimal model or $Z^{2}=-1$. Taking into account all possibilities
for the fiber $L$ to be a blow up of a reducible fiber of minimal
smooth elliptic fibration, we see that the equality $Z^{2}=-1$
holds, the curves $\hat{E}_{1}$ and $\hat{E}_{1}$ are sections of
the elliptic fibration $\psi$, but $\check{E}_{1}$ and
$\check{E}_{2}$ are contained in the fiber $L$. On the other hand,
the equivalences
$$
K_{Y}\sim_{\mathbb{Q}}\tau^{*}\big(\mathcal{O}_{\mathbb{P}(1,\,4,\,7,\,14)}(2)\vert_{S}\big)
-\frac{2}{7}\hat{E}_{1}-\frac{4}{7}\check{E}_{1}-\frac{2}{7}\hat{E}_{2}
-\frac{4}{7}\check{E}_{2}\sim_{\mathbb{Q}}2Z+E
$$
hold. Let $\gamma:Y\to\bar{Y}$ be the contraction of the curves
$Z$ and $E$. Then $\bar{Y}$ is smooth, the curve $\gamma(L)$ is a
fiber of type $III$ of the relatively minimal elliptic fibration
$\psi\circ\gamma^{-1}$, and the equivalence $K_{\bar{Y}}\sim 0$
holds. Therefore, the surface $S$ is birational to a $K3$ surface.

\emph{Case $\nn=83$.}

The surface $S$ is a hypersurface of degree $36$ in $
\mathbb{P}(1,4,11,18)\cong\mathrm{Proj}(\mathbb{C}[x_1,x_3,x_4,x_5]),
$ where $\mathrm{wt}(x_1)=1$, $\mathrm{wt}(x_3)=4$,
$\mathrm{wt}(x_4)=11$, and $\mathrm{wt}(x_5)=18$. Therefore, the
surface $S$ has a canonical singular point $Q$ of type
$\mathbb{A}_{1}$ given by the equations $x_1=x_4=0$ and an
isolated singular point $P$ at $(0:0:1:0)$.  The surface $S$ is
not quasismooth at the point $P$ which is not a rational singular
point of $S$, \emph{a posteriori}.

Let $\mathcal{P}$ be the pencil of curves on $S$ given by the
equations
$$
\lambda x_1^{4}+\mu x_3=0,
$$
where $(\lambda:\mu)\in\mathbb{P}^{1}$,  $C$ be a general curve in
$\mathcal{P}$, and $\nu:\hat{C}\to C$ be the normalization of the
curve $C$. Then the base locus of the pencil $\mathcal{P}$
consists of the point $P$ and $\mathcal{P}$ gives a rational map
$\phi:S\dasharrow\mathbb{P}^{1}$ whose general fiber is $C$. On
the other hand, the curve $C$ is a hypersurface of degree $36$ in
$\mathbb{P}(1,11,18)$. Therefore, the curve $\hat{C}$ is an
elliptic curve, and the birational map  $\nu$ is a bijection
because $C$ is a compactification of the affine curve
$$
C\cap \{x_1\ne 0\}\subset\mathbb{C}^{2},
$$ which is a double cover of $\mathbb{C}$
ramified at three points. In particular, we have
$\kappa(S)\leqslant 1$.

Let $\tau:Y\to S$ be the minimal resolution of singularities of
$S$. Then we have an elliptic fibration $\psi:Y\to\mathbb{P}^{1}$
such that $\psi=\phi\circ\tau$. We can identify a general fiber of
$\psi$ with the curve $\hat{C}$ and the normalization $\nu$ with
the restriction $\tau\vert_{\hat{C}}$. Therefore, there is exactly
one exceptional curve $Z$ of the resolution $\tau$ not contained
in a fiber of $\psi$. The curve $Z$ must be a section of $\psi$.

Let $F$ be the proper transform of the smooth rational curve in
the pencil $\mathcal{P}$ that is given by the equation $x_1=0$,
$E$ be the exceptional curve of the morphism $\tau$ that is mapped
to the point $Q$, and $E_{1},\ldots, E_{m}$ be the exceptional
curves of the birational morphism $\tau$ that are different from
the curves $Z$ and $E$. Then $\tau(E_{i})=\tau(Z)=P$ and the union
$$
F\cup E\cup E_{1}\cup\cdots\cup E_{m}
$$
lies in a single fiber $L$ of $\psi$. Moreover, the smooth
rational curve $F$ is a component of the fiber $L$ of multiplicity
$4$, the curve $E$ is rational, and $E^{2}=-2$. We have
$$
K_{Y}\sim_{\mathbb{Q}}2F+aE+\sum_{i=1}^{m}c_{i}E_{i},
$$
where $a$, $b$, $c_{i}$ are rational numbers. The elliptic
fibration $\psi$ is not relatively minimal, but the curve $F$ is
the only curve in the fiber $L$ whose self-intersection is $-1$.

Let $\xi:Y\to\bar{Y}$ be the birational morphism such that the
surface $\bar{Y}$ is the minimal model of the surface $Y$ and
$\eta=\psi\circ\xi^{-1}$. Then $\eta:\bar{Y}\to\PP^1$ is a
relatively minimal elliptic fibration.

Let $\bar{L}=\xi(L)$. Then
$K_{\bar{Y}}\sim_{\mathbb{Q}}\gamma\bar{L}$ for some rational
number $\gamma\geqslant 0$. Hence, we have
$$
K_{Y}\sim_{\mathbb{Q}}\xi^{*}(\gamma\bar{L})+\alpha F+\beta E+\sum_{i=1}^{m}\delta_{i}E_{i},%
$$
where $\alpha$, $\beta$, $\delta_{i}$ are non-negative integer
numbers. Because the birational morphism $\xi$ must contract the
curves $F$ and $E$, we see that $\alpha\geqslant 2$,
$\beta\geqslant 1$. Also, the inequality $\delta_{i}\ne 0$ holds
if and only if the curve $E_{i}$ is contracted by $\xi$. Moreover,
the equality $\alpha=2$ implies that the only curves contracted by
$\xi$ are $F$ and $E$. Hence, the inequality $\gamma\geqslant 0$
and the equivalence
$$
2F+aE+\sum_{i=1}^{m}c_{i}E_{i}\sim_{\mathbb{Q}}\xi^{*}(\gamma\bar{L})+
\alpha F+\beta E+\sum_{i=1}^{m}\delta_{i}E_{i}%
$$
imply that $\gamma=0$, $\alpha=2$, and $m>0$. In particular, the
surface $\bar{Y}$ is either a $K3$ surface or an Enriques surface.
On the other hand, the only possible multiple fiber of the
elliptic fibration $\eta$ is the fiber $\bar{L}$, which implies
that $\bar{Y}$ is a $K3$ surface.

Therefore, we have proved
Theorem~\ref{proposition:existence-of-K3-fibration}. In addition,
we have  shown that $X$ is birational to a fibration whose general
fiber is an elliptic $K3$ surface if $\nn\not\in\{3, 60, 75, 87,
93\}$.

We conclude the section with one remark.
\begin{remark}{\normalfont In the proof of Case $\nn=83$ the equality $\alpha=2$ and the
fact that $F$ is a component of $L$ of multiplicity $4$ imply that
$\bar{L}$ is an elliptic fiber of type $I_{r}^{*}$, while the
birational morphism $\xi$ is the composition of the blow up at a
point of the component of the fiber $\bar{L}$ of multiplicity $2$
and the blow up at the intersection point of the proper transform
of the component of multiplicity $2$ with the exceptional curve on
the first blow up.

It was pointed out to us by D.\,Stepanov that one can explicitly
resolve the singularity of the surface $S$ at the point $P$ to
prove that the surface $S$ is birationally equivalent to a smooth
$K3$ surface. Indeed, the surface $S$ can be locally given near
$P$ by the equation
$$
x^2+y^3+z^9=0\subset\mathbb{C}^{3}\slash\,\mathbb{Z}_{11}(7,4,1)\,%
$$
where $P=(0,0,0)$.

Let $\sigma_{1}$ be the weighted blow up of
$\mathbb{C}^{3}/\,\mathbb{Z}_{11}(7,4,1)$ at the singular point
$P$ with the weight $\frac{1}{11}(10,3,1)$. Then the blown up
variety is covered by $3$ affine charts, the first chart is
isomorphic to $\mathbb{C}^3/\mathbb{Z}_{10}(1,-3,-1)$, and in the
first chart $\sigma_{1}$ is given by
$$
x=x^{10/11},\ y=x^{3/11}y,\ z=x^{1/11}z,
$$
where we denote the coordinates on
$\mathbb{C}^3/\mathbb{Z}_{10}(1,-3,-1)$ by the same letters $x$,
$y$, $z$ as the coordinates on
$\mathbb{C}^{3}/\,\mathbb{Z}_{11}(7,4,1)$. The full transform of
$S$ is given by the equation
$$
x^{20/11}+x^{9/11}y^3+x^{9/11}z^9=0,
$$
but the strict transform $\bar{S}$ of the surface $S$ is given by
the equation
$$
x+y^3+z^9=0\subset\mathbb{C}^3/\mathbb{Z}_{10}(1,-3,-1),
$$
and the exceptional divisor
$$
x=0=y^3+z^9=\prod_{i=1}^{3}(y+\varepsilon^{i} z^9)
$$
consists of $3$ smooth rational curves $\bar{E}_{1}$,
$\bar{E}_{2}$, $\bar{E}_{3}$ that intersect at the singular point
$(0,0,0)$, where $\varepsilon$ is a primitive cubic root of unity.
Moreover, the surface $\bar{S}$ has quotient singularity of type
$\frac{1}{10}(-3,-1)$ at the singular point $(0,0,0)$.

In the second chart that is  isomorphic to
$\mathbb{C}^3/\mathbb{Z}_{3}(-1,2,-1)$, the strict transform of
$S$ is given by the equation $x^2y+1+z^9=0$, and in the third
chart that is  isomorphic to $\mathbb{C}^3$, the strict transform
of the surface $S$ is given by $x^2z+y^3+1=0$, which imply that
they are nonsingular.

We have a surface $\bar{S}$ that is locally isomorphic to
$\mathbb{C}^2\slash\,\mathbb{Z}_{10}(-3,-1)$ and we have $3$
smooth rational curves on $\bar{S}$ given by the equation
$$
\prod_{i=1}^{3}(x+\varepsilon^i y^3)=0,
$$
where $x$ and $y$ are local coordinates on
$\mathbb{C}^2\slash\,\mathbb{Z}_{10}(-3,-1)$.

Let $\sigma_{2}$ be the weighted blow up of the surface $\bar{S}$
at the point $(0,0)$ with the weight $\frac{1}{10}(1,7)$. The
blown up variety is covered by $2$ charts.  The first chart is
$\mathbb{C}^2$ and it does not contain the strict transforms of
the curve $\bar{E_{i}}$. The second chart is isomorphic to
$\mathbb{C}^2\slash\,\mathbb{Z}_{7}(-1,3)$ and in this chart the
weighted blow up $\sigma_2$ is given by the formulas
$x=y^{1/10}x$, $y=y^{7/10}$ but the strict transform of the curve
$\bar{E_{i}}$ is given by the equation $x+\varepsilon^i y^2=0$,
where the exceptional divisor $\bar{Z}$ of the weighted blow up
$\sigma_2$ is given by $y=0$.

Now let $\sigma_3$ be the weighted blow up at the origin of the
last considered chart with the weight $\frac{1}{7}(2,1)$. In the
first chart $\mathbb{C}^2\slash\,\mathbb{Z}_{2}(1,1)$, the
equation of the proper transform of the curve $\bar{E_{i}}$ is
$1+\varepsilon^i y^2=0$, the equation of the proper transform of
$\bar{Z}$ is $y=0$, and the exceptional divisor $\bar{E_4}$ of
$\sigma_3$ is given by $x=0$, but the second chart of $\sigma_3$
is nonsingular.

Let $\sigma_{4}$ be the blow up of
$\mathbb{C}^2\slash\,\mathbb{Z}_{2}(1,1)$ with the weight
$\frac{1}{2}(1,1)$ and let $\bar{E}_{5}$ be the exceptional
divisor of $\sigma_{4}$. Then $\sigma_{4}$ resolves the
singularity of $S$ in a neighborhood of the point $P$ and after
blowing up the point $Q$ of $S$ we get our minimal resolution
$\tau:Y\to S$.

Let $E_{i}$ and $Z$ be the proper transforms of the irreducible
curves $\bar{E}_{i}$ and $\bar{Z}$ on the nonsingular surface $Y$,
respectively. Then $E_{4}^{2}=-4$, $Z^{2}=E_{i\ne
4}^{2}=E^{2}=-2$, where
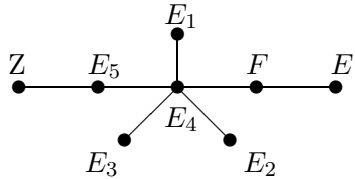
\begin{center}
\begin{picture}(360,80)(0,0)
\multiput(100,40)(30,0){5}{\circle*{5}}\put(140,20){\circle*{5}}
\put(160,60){\circle*{5}}\put(180,20){\circle*{5}}
\put(100,40){\line(1,0){120}}\put(160,60){\line(0,-1){20}}\put(160,40){\line(1,-1){20}}
\put(160,40){\line(-1,-1){20}} \put(96,45){Z}
\put(126,45){$E_{5}$}\put(155,64){$E_{1}$} \put(155,25){$E_{4}$}
\put(186,45){$F$} \put(218,45){$E$} \put(185,8){$E_2$}
\put(125,8){$E_{3}$}
\end{picture}
\end{center}
is the dual graph of the rational curves $Z, E_{1},\ldots, E_{5}$,
$F$, and $E$. In particular, the  fiber $\bar{L}$ is of type
$I^*_0$.}\end{remark}

\section{Birational automorphisms.}
\label{section:birational-automorphisms}

The group $\mathrm{Bir}(X)$ of birational automorphisms is
generated by biregular automorphisms and a finite set of
birational involutions $\tau_{1},\ldots,\tau_{\ell}$ that are
described in \cite{CPR}. To be precise, we have an exact sequence
of groups
$$
1\to \Gamma_X\to \mathrm{Bir}(X)\to \mathrm{Aut}(X)\to 1,
$$
where the group $\Gamma_X$ is a subgroup of $\mathrm{Bir}(X)$
generated by a finite set of distinct birational involutions
$\tau_{1},\ldots,\tau_{\ell}$.

In this section we describe the group $\Gamma_X$ with group
presentations. When the number $\ell$ of generators of $\Gamma_X$
is $0$, namely, the group $\Gamma_X$ is trivial, $\Bir
(X)=\Aut(X)$, and hence the $3$-fold $X$ is birationally
superrigid.  When the number $\ell$ of generators of $\Gamma_X$ is
$1$, the group $\Gamma_X$ is the group of order $2$, \emph{i.e.},
$\ZZ/2\ZZ$. Therefore, we may assume that $\ell\geqslant 2$ to
prove Theorem~\ref{theorem:birational-autmorphisms}. Throughout
this section, a relation of involutions means one different from
the trivial relation, \emph{i.e.}, $\tau_i^2=1$.

First of all, we present the following important observation:

\begin{lemma}\label{lemma:One-point}Suppose that the set $\mathbb{CS}(X,
\lambda\mathcal{M})$ contains at most one element, where
$\mathcal{M}$ is a linear system without fixed components on $X$
and $\lambda$ is a positive rational number such that the divisor
$-(K_X+\mathcal{M})$ is ample. Then there is no relation among
$\tau_1,\cdots, \tau_\ell$.
\end{lemma}

\begin{proof}
See Proposition~2.2 and Lemma~2.3 in \cite{Pu88b}. They show the
condition implies a given birational automorphism is
untwisted\footnote{
Fix a very ample linear system $\mathcal{H}$ on $X$. Let
$\phi:X\dasharrow X$ be a birational automorphism such that
$\phi^{-1}(\mathcal{H})\subset |-rK_X|$. We say that an involution
$\tau$ of $X$ untwists the map $\phi$ if $(\phi\circ\tau
)^{-1}(\mathcal{H})\subset |-r'K_X|$ for some $r'<r$. More
generally, for a log pair $(X,\frac{1}{r}\mathcal{M})$ with
$\mathcal{M}\sim_\QQ-K_X$ that is not terminal we also say that an
involution $\tau$ of $X$ untwists a maximal singularity of
$(X,\frac{1}{r}\mathcal{M})$ if $\tau(\mathcal{M})\sim_\QQ-r'K_X$
for some  $r'<r$. For more generalized detail, refer to
\cite{CPR}.} by the involutions $\tau_1, \ldots,\tau_\ell$ in a
unique way (See also \cite{CPR}).
\end{proof}

Note that the assumption $\ell\geq 2$ implies that
$$
\nn\in\{4, 5, 6, 7, 9, 12, 13, 15, 17, 20, 23, 25, 27, 30, 31, 33,
36, 38, 40, 41, 42, 44, 58, 61, 68, 76\}.
$$

\begin{lemma}
\label{lemma:free-group1} Suppose that $\nn\in\{6, 15, 23, 30, 36,
40, 41, 42, 44, 61, 68, 76\}$. Then $\Gamma_X$ is the free product
of two involutions $\tau_1$ and $\tau_2$.
\end{lemma}

\begin{proof}
Suppose that $\nn=36$. Then the hypersurface$X$ is a sufficiently
general hypersurface in $\mathbb{P}(1,1,4,6,7)$ of degree $18$
with $-K_{X}^{3}=3/28$. It has three singular points, namely, the
point $P_{1}$ that is a quotient singularity of type
$\frac{1}{2}(1,1,1)$, the point $P_{2}$ that is a quotient
singularity of type $\frac{1}{4}(1,1,3)$, and the point $P_{3}$
that is a quotient singularity of type $\frac{1}{7}(1,1,6)$.

Suppose that the group $\Gamma_X$ is not the free product of the
involutions $\tau_{1}$ and $\tau_{2}$. Then there is a linear
system $\mathcal{M}$ without fixed components on the hypersurface
$X$ such that  the set $\mathbb{CS}(X, \lambda\mathcal{M})$
contains at least two subvariety of the hypersurface $X$, where
$\lambda$ is a positive rational number such that the divisor
$-(K_{X}+\lambda\mathcal{M})$ is ample. Therefore, it follows from
\cite{CPR} that $\mathbb{CS}(X,
\lambda\mathcal{M})=\{P_{2},P_{3}\}$.

The hypersurface $X$ can be given by the quasihomogeneous equation
of degree $18$
$$
x_3^{3}x_4+x_3^{2}g(x_1,x_2,x_4,x_5))+x_3h(x_1,x_2,x_4,x_5)+q(x_1,x_2,x_4,x_5)=0\subset\PP(1,1,4,6,7),
$$
where $f$, $g$, $h$, and $q$ are quasihomogeneous polynomials.
Then the point $P_{2}$ is  located at $(0:0:1:0:0)$ and the point
$P_{3}$ at $(0:0:0:0:1)$.

Let $\xi:X\dasharrow\mathbb{P}^{7}$ be the rational map that is
given by the linear subsystem of the linear system $|-6K_{X}|$
consisting of the  divisors
$$
\mu x_4+\sum_{i=0}^{6}\lambda_{i}x^{i}_1x_2^{6-i}=0,
$$
where
$(\mu:\lambda_{0}:\lambda_{1}:\lambda_{2}:\lambda_{3}:\lambda_{4}:\lambda_{5}:\lambda_{6})\in\mathbb{P}^{7}$.
Then the rational map $\xi$ is not defined at the points $P_{2}$
and $P_{3}$, the closure of the image of the rational map $\xi$ is
the surface $\mathbb{P}(1,1,6)$, and a general fiber of the map
$\xi$ is an elliptic curve. There is a commutative diagram
$$
\xymatrix{
&&W\ar@{->}[dl]_{\beta_{3}}\ar@{->}[dr]^{\beta_{2}}\ar@{->}[rrrrrrrd]^{\omega}&&&&&\\
&U_{2}\ar@{->}[dr]_{\alpha_{2}}&&U_{3}\ar@{->}[dl]^{\alpha_{3}}&&&&&&\mathbb{P}(1,1,6),\\
&&X\ar@{-->}[rrrrrrru]_{\xi}&&&&&&&}
$$
where $\alpha_{2}$ is the Kawamata blow up at the singular point
$P_{2}$, $\alpha_{3}$ is the Kawamata blow up at the point
$P_{3}$, $\beta_{2}$ is the Kawamata blow up at the point
$\alpha_3^{-1}(P_{2})$, $\beta_{3}$ is the Kawamata blow up at the
point  $\alpha_2^{-1}(P_{3})$, and $\omega$ is an elliptic
fibration.

Let $S$ be the proper transform on the $3$-fold $W$ of a general
surface of the linear system $\mathcal{M}$ and $C$ be a general
fiber of the fibration $\omega$. The inequality $S\cdot C<0$
follows from \cite{Ka96}. However, it is a contradiction because
$\omega$ is an elliptic fibration.

Suppose that $\nn=44$. Then $X$ is a general hypersurface in
$\mathbb{P}(1,2,5,6,7)$ of degree $20$ with $-K_{X}^{3}=1/21$. The
singularities of the hypersurface $X$ consist of the points
$P_{1}$, $P_{2}$, $P_{3}$ that are quotient singularities of type
$\frac{1}{2}(1,1,1)$, the point $P_{4}$ that is a quotient
singularity of type $\frac{1}{6}(1,1,5)$, and the point $P_{5}$
that is a quotient singularity of type $\frac{1}{7}(1,2,5)$.
Moreover, there is a commutative diagram
$$
\xymatrix{
&&Y\ar@{->}[dl]_{\beta_{5}}\ar@{->}[dr]^{\beta_{4}}\ar@{->}[rrrrrrrd]^{\eta}&&&&&\\
&U_{4}\ar@{->}[dr]_{\alpha_{4}}&&U_{5}\ar@{->}[dl]^{\alpha_{5}}&&&&&&\mathbb{P}(1,2,5),\\
&&X\ar@{-->}[rrrrrrru]_{\psi}&&&&&&&}
$$
where $\psi$ is the natural projection, $\alpha_{4}$ is the
weighted blow up at the singular point $P_{4}$ with weights
$(1,1,5)$, $\alpha_{5}$ is the weighted blow up at the point
$P_{5}$ with weights $(1,2,5)$, $\beta_{4}$ is the weighted blow
up with weights $(1,1,5)$ at the point $\alpha_5^{-1}(P_4)$,
$\beta_{5}$ is the weighted blow up with weights $(1,2,5)$ at the
point $\alpha^{-1}_4(P_5)$, and $\eta$ is an elliptic fibration.
It follows from \cite{CPR} that
$$
\mathbb{CS}(X, \lambda\mathcal{M})=\{P_{4},P_{5}\},
$$
and we can proceed as in the previous case to obtain a
contradiction.

In the case when $\nn\in\{6, 15, 23, 30, 40, 41, 42, 61, 68, 76\}$
we can obtain a contradiction in the same way as in the case
$\nn=44$.
\end{proof}

\begin{lemma}
\label{lemma:free-group-4-9-17-27} Suppose that $\nn\in\{4, 9, 17,
27\}$. Then
$\tau_{1}\circ\tau_{2}\circ\tau_{3}=\tau_{3}\circ\tau_{2}\circ\tau_{1}$
is the only relation among the birational involutions
$\tau_{1},\tau_{2}$, and $\tau_{3}$.
\end{lemma}

\begin{proof} It follows from \cite{CPR} that $\ell=3$, $a_{4}=a_{5}$, and $d=3a_{4}$.
A general fiber of the projection
$\psi:X\dasharrow\mathbb{P}(1,a_{2},a_{3})$ is a smooth elliptic
curve. Moreover, the hypersurface $X$ has singular points $P_{1}$,
$P_{2}$, $P_{3}$ of index $a_{4}$ which are the points of the
indeterminacy of the map $\psi$.

Let $\pi:V\to X$ be the Kawamata blow up at the points $P_{1}$,
$P_{2}$, $P_{3}$. We also  let $E_{i}$ be the exceptional divisor
of $\pi$ dominating $P_{i}$ and $\phi=\psi\circ\pi$. Then $\pi$ is
a resolution of indeterminacy of the rational map $\psi$, the
divisors $E_{1}$, $E_{2}$, $E_{3}$ are sections of $\phi$, the
equivalence
$$
-K_{V}\sim_{\mathbb{Q}}\pi^{*}(-K_{X})-\frac{1}{a_{4}}E_{1}-\frac{1}{a_{4}}E_{2}-\frac{1}{a_{4}}E_{3}
$$
holds, the linear system $|-a_{3}a_{4}a_{5}K_{V}|$ is free and
lies in the fibers of $\phi$.

Let $\mathbb{F}$ be the field of rational functions on
$\mathbb{P}(1,a_{2},a_{3})$ and $C$ be a generic fiber of the
elliptic fibration $\phi$ considered as an elliptic curve over
$\mathbb{F}$. Then the section $E_{j}$ of the elliptic fibration
$\phi$ can be considered as an $\mathbb{F}$-rational point of the
elliptic curve $C$.

One can show using Lemma~\ref{lemma:elliptic-surfaces} that
$\mathbb{F}$-rational points $E_{1}$, $E_{2}$, $E_{3}$ are
$\mathbb{Z}$-linearly independent in the group $\mathrm{Pic}(C)$,
but we never use the latter in the rest of the proof.

By our construction, the curve $C$ is a hypersurface of degree
$3a_{4}$ in $\mathbb{P}(1,a_{4},a_{4})\cong\mathbb{P}^{2}$, which
implies that the curve $C$ can be naturally identified with a
cubic curve in $\mathbb{P}^{2}$ such that the points $E_{1}$,
$E_{2}$, $E_{3}$ lie on a single line in $\mathbb{P}^{2}$.

Let $\sigma_{i}$ be the involution of the curve $C$ that
interchanges the fibers of the projection of the curve $C$ from
the point $E_{i}$. Then $\sigma_{i}$ can be also considered as a
birational involution of the 3-fold $V$ such that
$$
\sigma_{i}=\pi^{-1}\circ\tau_{i}\circ\pi\in\mathrm{Bir}(V).
$$

Consider the curve $C$ as a group scheme. Let $Q_{k}$ be  the
point $(E_{i}+E_{j})\slash 2$ on the elliptic curve $C$, where
$\{i,j\}=\{1,2,3\}\setminus \{k\}$. Then the involution
$\sigma_{k}$ is the reflection of the elliptic curve $C$ at the
point $Q_{k}$ because the points $E_{1}$, $E_{2}$, $E_{3}$ are
$\mathbb{Z}$-linearly independent, which implies that $Q_{1}$,
$Q_{2}$, $Q_{3}$ are $\mathbb{Z}$-linearly independent and the
compositions
$$
\sigma_{2}\circ\sigma_{1}\circ\sigma_{3},\
\sigma_{1}\circ\sigma_{2}\circ\sigma_{3},\
\sigma_{1}\circ\sigma_{3}\circ\sigma_{2}
$$
are reflections at $E_{1}$, $E_{2}$, $E_{3}$ respectively. Thus,
we have the identity
$$
\tau_{1}\circ\tau_{2}\circ\tau_{3}=\tau_{3}\circ\tau_{2}\circ\tau_{1},
$$
which implies the similar identities that can be obtained from
$\tau_{1}\circ\tau_{2}\circ\tau_{3}=\tau_{3}\circ\tau_{2}\circ\tau_{1}$
by a permutation of the elements in the set $\{1,2,3\}$.

It follows from \cite{CPR} that for any linear system
$\mathcal{M}$ on the hypersurface $X$ having no fixed components,
the singularities of the log pair $(X, \frac{1}{r} \mathcal{M})$
are canonical in the outside of the points $P_{1}$, $P_{2}$,
$P_{3}$, where $r$ is the natural number such that
$\mathcal{M}\sim_{\mathbb{Q}} -rK_{X}$. Moreover, when the
singularities of the log pair $(X, \frac{1}{r} \mathcal{M})$ are
not canonical at the point $P_{i}$, we have
$$
\frac{1}{r}\mathcal{B}\sim_{\mathbb{Q}}\pi^{*}(\frac{1}{r}\mathcal{M})
-m_{1}E_{1}-m_{2}E_{2}-m_{3}E_{3},
$$
where $\mathcal{B}$ is the proper transform of $\mathcal{M}$ on
$V$ and $m_{i}>1/a_{4}$. We have the inequality
$$
m_{1}+m_{2}+m_{3}\leqslant \frac{3}{a_{4}},
$$
which implies that the linear system $\mathcal{B}$ lies in the
fibers of the elliptic fibration $\phi$ if the equality
$m_{1}+m_{2}+m_{3}=3/a_{4}$ holds.

When the inequality $m_{i}>1/a_{4}$ holds, the birational
involution $\tau_{i}$ untwists the maximal singularity of the log
pair $(X, \frac{1}{r} \mathcal{M})$ at the point $P_{i}$, namely,
the equivalence
$$
\tau_{i}(\mathcal{M})\sim_{\mathbb{Q}}-r^{\prime}K_{X}
$$
holds for some natural number $r^{\prime}<r$. Similarly, the
involution $\tau_{i}\circ\tau_{k}\circ\tau_{j}$ untwists the
maximal singularities of the log pair $(X, \frac{1}{r}
\mathcal{M})$ at the points $P_{i}$ and $P_{j}$ simultaneously
when the inequalities $m_{i}>1/a_{4}$ and $m_{j}>1/a_{4}$ holds
for $i\ne j$, where $k\in\{1,2,3\}\setminus\{i,j\}$.

Now we can use the arguments of the proof of Theorem~7.8 in
Section V of \cite{Ma72} to prove that the identity
$\tau_{1}\circ\tau_{2}\circ\tau_{3}=\tau_{3}\circ\tau_{2}\circ\tau_{1}$
is the only relation    among our birational involutions
$\tau_{1}$, $\tau_{2}$, and $\tau_{3}$. However, it should be
pointed out that the arguments of the proof of Theorem~7.8 in
Section V of \cite{Ma72} are too sophisticated for our
purposes\footnote{The following arguments are due to A.Borisov.
Let $W$ be a composition of $\sigma_{1}$, $\sigma_{2}$,
$\sigma_{3}$ such that $W$ is the identity map of the elliptic
curve $C$ and $W$ does not contain squares of $\sigma_{i}$. Then
we can show that $W$ has even number of entries and each entry
appears the same number of times in the even and odd position, and
we can use the identity
$\sigma_{1}\circ\sigma_{2}\circ\sigma_{3}=\sigma_{3}\circ\sigma_{2}\circ\sigma_{1}$
to make $\sigma_{3}$ jump $2$ spots left or right. Shifting the
last $\sigma_{3}$ in the odd position in $W$ that is followed not
right away by $\sigma_{3}$ in the even position, we can collapse
them and get a composition of $\sigma_{1}$, $\sigma_{2}$,
$\sigma_{3}$ having a smaller number of entries. Therefore, the
only relation among the involutions $\sigma_{1}$, $\sigma_{2}$,
$\sigma_{3}$ is the identity
$\sigma_{1}\circ\sigma_{2}\circ\sigma_{3}=\sigma_{3}\circ\sigma_{2}\circ\sigma_{1}$.}.
\end{proof}

\begin{lemma}
\label{lemma:free-group-7} Suppose that $\nn=7$. Then there are no
relations among $\tau_{1},\ldots,\tau_{5}$.
\end{lemma}

\begin{proof}
The 3-fold $X$ is a general hypersurface in
$\mathbb{P}(1,1,2,2,3)$ of degree $8$ which has singular points
$P_{1},\ldots,P_{4}$ of type ${\frac{1}{2}}(1,1,1)$ and a singular
point $Q$ of type ${\frac{1}{3}}(1,2,1)$.

Let $\alpha_{i}:V_{i}\to X$ be the weighted blow up of $X$ at the
singular points $P_{i}$ and $Q$ with the weights
$\frac{1}{2}(1,1,1)$ and $\frac{1}{3}(1,2,1)$, respectively. Then
$$
K_{V_{i}}\sim_{\mathbb{Q}}\alpha_{i}^{*}(K_{X})+{\frac{1}{2}}E_{i}+{\frac{1}{3}}F_{i},%
$$
where $E_{i}$ and $F_{i}$ are the exceptional divisors of the
birational morphism $\alpha_{i}$ dominating the singular points
$P_{i}$ and $Q$, respectively. The linear system $|-2K_{V_{i}}|$
induces the morphism
$$
\psi_{i}:V_{i}\to\mathbb{P}(1,1,2),
$$
which is an elliptic fibration. Moreover, the divisor $E_{i}$ is a
$2$-section of the fibration $\psi_{i}$, while the divisor $F_{i}$
is a section of $\psi_{i}$. Up to relabelling, the birational
involutions $\tau_{1},\ldots,\tau_{5}$ can be constructed as
follows: the involution $\tau_{i}$ is induced by the reflection of
a general fiber of the morphism $\psi_{i}$ at the section $F_{i}$
but the involution $\tau_{5}$ is induced by the natural projection
$X\dasharrow\mathbb{P}(1,1,2,2)$.

Let $\mathcal{M}$ be a linear system on $X$ without fixed
components such that $\mathcal{M}\sim_{\mathbb{Q}} -rK_{X}$ for
some natural number $r$. Then the singularities of the log pair
$(X, \frac{1}{r} \mathcal{M})$ are canonical in the outside of the
points $P_{1},\ldots,P_{4},Q$ due to \cite{CPR}.

Let $\mathcal{B}_{i}$ be the proper transform of $\mathcal{M}$ on
$V_{i}$. Then
$$
\mathcal{B}_i\sim_{\mathbb{Q}}\alpha_{i}^{*}(\mathcal{M})-m_{i}E_{i}-mF_{i},
$$
where $m_{i}$ and $m$ are positive rational numbers. Moreover, the
log pair $(X, \frac{1}{r} \mathcal{M})$ is not canonical at the
point $P_{i}$ if and only if $m_{i}>r/2$. On the other hand, the
singularities of the log pair $(X, \frac{1}{r} \mathcal{M})$ are
not canonical at $Q$ if and only if $m>r/3$. Now intersecting the
linear system $\mathcal{B}_{i}$ with a sufficiently general fiber
of $\psi_{i}$, we see that
$$
2m_{i}+m\leqslant \frac{4r}{3}.
$$

The equivalence
$\tau_{i}(\mathcal{M})\sim_{\mathbb{Q}}-r^{\prime}K_{X}$ holds for
some natural number $r^{\prime}$. Moreover, the inequality
$r^{\prime}<r$ holds if $m_{i}>r/2$ when $i=1,\ldots, 4$ or if
$m>r/3$ when $i=5$, namely, the involutions
$\tau_{1},\ldots,\tau_{5}$ untwist the maximal singularities of
the log pair $(X, \frac{1}{r} \mathcal{M})$.

In order to prove that the involutions $\tau_{1},\ldots,\tau_{5}$
do not have any relation, it is enough to prove that  the
singularities of the log pair $(X, \frac{1}{r} \mathcal{M})$ are
not canonical at at most one point by Lemma~\ref{lemma:One-point}.
However, the inequality $2m_{i}+m\leqslant 4r/3$ implies that the
log pair $(X, \frac{1}{r} \mathcal{M})$ is canonical at one of the
singular points $P_{i}$ and $Q$. To conclude the proof, therefore,
we must show that for $i\ne j$ the singularities of the log pair
$(X, \frac{1}{r} \mathcal{M})$ are canonical at one of the points
$P_{i}$ and $P_{j}$.

Suppose that the log pair $(X, \frac{1}{r} \mathcal{M})$ is not
canonical at the points  $P_{1}$ and $P_{2}$.  Let $S$ be a
general surface in $|-K_{X}|$ and $C$ be the base curve of
$|-K_{X}|$. Then $S$ is a $K3$ surface whose singular points are
the singular points of $X$. Moreover, the point $P_{i}$ is a
singular point of type $\mathbb{A}_{1}$ on the surface $S$ and the
point $Q$ is a singular point of type $\mathbb{A}_{2}$ on $S$.

The curve $C$ is a smooth curve passing through the points
$P_{1},\ldots,P_{4}$, and $Q$. We have
$$
\mathcal{M}\vert_{S}=\mathcal{P}+\mathrm{mult}_{C}(\mathcal{M})C,
$$
where $\mathcal{P}$ is a linear system without fixed components.
The inequality $\mathrm{mult}_{C}(\mathcal{M})\leqslant r$ holds;
otherwise the log pair $(X, \frac{1}{r} \mathcal{M})$ would not be
canonical at the point $Q$ by \cite{Ka96}.

Let $\pi:Y\to S$ be the composition of blow ups of the singular
points $P_{1}$ and $P_{2}$, $G_{i}$ be the exceptional divisor of
$\pi$ dominating $P_{i}$, $\bar{C}$ be the proper transform on the
surface $Y$ of the curve $C$, and $\mathcal{H}$ be the proper
transform on the surface $Y$ of the linear system $\mathcal{P}$.
Then
$$
\mathcal{H}+\mathrm{mult}_{C}(\mathcal{M})\bar{C}
\sim_{\mathbb{Q}}\pi^{*}\big(-rK_{X}\vert_{S}\big)-\bar{m}_{1}G_{1}-\bar{m}_{2}G_{2},
$$
where $\bar{m}_{i}\geqslant m_{i}>r/2$. However, we have
$\bar{C}^{2}=-1/3$ on the surface $Y$ and we see that
$$
-\frac{r}{3}\leqslant\big(\mathcal{H}+\mathrm{mult}_{C}(\mathcal{M})
\bar{C}\big)\cdot\bar{C}\leqslant \frac{2r}{3}-\bar{m}_{1}-\bar{m}_{2}<-\frac{r}{3},%
$$
which is a contradiction.
\end{proof}

\begin{lemma}
\label{lemma:free-group-20} Suppose that $\nn=20$. Then there are
no relations among $\tau_{1},\tau_{2},\tau_{3}$.
\end{lemma}

\begin{proof}
We have a general hypersurface $X\subset\mathbb{P}(1,1,3,4,5)$
given by
$$
x_3^{4}f_{1}(x_1,x_2)+x_3^{3}f_{4}(x_1,x_2,x_4)+x_3^{2}f_{7}(x_1,x_2,x_4,x_5)+$$$$
+x_3f_{10}(x_1,x_2,x_4,x_5)+f_{13}(x_1,x_2,x_4,x_5)=0,
$$
where  $f_{i}$ is a general quasihomogeneous polynomial of degree
$i$. The 3-fold $X$ has $3$ singular points at $P=(0:0:1:0:0)$,
$Q=(0:0:0:1:0)$, $O=(0:0:0:0:1)$ and a general fiber of the
natural projection of $X$ to $\mathbb{P}(1,1,3)$ is an elliptic
curve. However, a general fiber of the natural projection of $X$
to $\mathbb{P}(1,1,4)$ may not be an elliptic curve.

Let us take $t=x_3f_{1}(x_1,x_2)+f_{4}(x_1,x_2,x_4)$ as a
homogeneous variable of weight $4$ instead of the homogeneous
variable $x_4$. Then the hypersurface $X$ is given by the equation
$$
x_3^{3}t+x_3^{2}g_{7}(x_1,x_2,t,x_5)+zg_{10}(x_1,x_2,t,x_5)+g_{13}(x_1,x_2,t,x_5)=0,
$$
where  $g_{i}$ is a sufficiently general quasihomogeneous
polynomial of degree $i$. A general fiber of the natural
projection of $X$ to $\mathbb{P}(1,1,4)$ is an elliptic curve.

Up to relabelling, the involutions $\tau_{1}$, $\tau_{2}$,
$\tau_{3}$ can be constructed as follows:
\begin{itemize}
\item the birational involution $\tau_{1}$ is induced by the
reflection of a general fiber of the natural projection
$X\dasharrow\mathbb{P}(1,1,4)$ at the point $O$; %
\item the birational involution $\tau_{2}$ is induced by the
reflection of a
general fiber of the natural projection $X\dasharrow\mathbb{P}(1,1,3)$ at the point $O$; %
\item the birational involution $\tau_{3}$ is induced by the
reflection of a general fiber of the natural projection
$X\dasharrow\mathbb{P}(1,1,3)$ at the point $Q$ but the involution
$\tau_{3}$ is also induced by the natural
projection $X\dasharrow\mathbb{P}(1,1,3,4)$.%
\end{itemize}

Let $\mathcal{M}$ be any linear system on $X$ without fixed
components such that $\mathcal{M}\sim_{\mathbb{Q}} -rK_{X}$ for
some natural number $r$. Then the singularities of the log pair
$(X, \frac{1}{r} \mathcal{M})$ are canonical in the outside of the
points $P$, $Q$, $O$ due to \cite{CPR}, and the equivalence
$\tau_{i}(\mathcal{M})\sim_{\mathbb{Q}}-r^{\prime}K_{X}$ holds for
some natural number $r^{\prime}<r$ in the following cases:
\begin{itemize}
\item the log pair $(X, \frac{1}{r}
\mathcal{M})$ are not canonical at the point $P$ and $i=1$;%
\item the log pair $(X, \frac{1}{r}
\mathcal{M})$ are not canonical at the point $Q$ and $i=2$; %
\item the log pair $(X, \frac{1}{r} \mathcal{M})$ are not
canonical at the point $O$ and $i=3$.
\end{itemize}

In order to prove that the involutions $\tau_{1}$, $\tau_{2}$,
$\tau_{3}$ are not related by any relation, by
Lemma~\ref{lemma:One-point} it is enough to show that the
singularities of $(X, \frac{1}{r} \mathcal{M})$ are not canonical
at at most one point.

Suppose that $(X, \frac{1}{r} \mathcal{M})$ is not canonical at
the points  $P$ and $O$. Let $\alpha:V\to X$ be the Kawamata blow
up at the points $P$ and $O$. Then
$$
K_{V}\sim_{\mathbb{Q}}\alpha^{*}(K_{X})+{\frac{1}{3}}E+{\frac{1}{5}}F,%
$$
where $E$ and $F$ are the exceptional divisors of the birational
morphism $\alpha$ dominating the singular points $P$ and $O$,
respectively. The linear system $|-4K_{V}|$ does not have base
points and induces the morphism $\psi:V\to\mathbb{P}(1,1,4)$ which
is an elliptic fibration. The divisor $F$ is a section of $\psi$
and the divisor $E$ is a $2$-section of $\psi$. Let $\mathcal{B}$
be the proper transform of the linear system $\mathcal{M}$ on the
3-fold $V$. Then
$$
\mathcal{B}\sim_{\mathbb{Q}}\alpha^{*}(\mathcal{M})-aE-bF,
$$
where $a$ and $b$ are rational numbers such that $a>r/3$ and
$b>r/5$. Intersecting the linear system $\mathcal{B}$ with a
sufficiently general fiber of $\psi$, we see that
$$
2a+b\leqslant\frac{52r}{60},
$$
which is impossible because $a>r/3$ and $b>r/5$.

We next suppose that the singularities of the log pair $(X,
\frac{1}{r} \mathcal{M})$ are not canonical at the singular points
$Q$ and $O$. Let $\gamma:W\to X$ be the composition of the
weighted blow ups at the points $Q$ and $O$ with weights $(1,1,3)$
and $(1,1,4)$, respectively. Then
$$
K_{W}\sim_{\mathbb{Q}}\gamma^{*}(K_{X})+{\frac{1}{4}}G+{\frac{1}{5}}H,%
$$
where $G$ and $H$ are the $\gamma$-exceptional divisors dominating
the singular points $Q$ and $O$, respectively. Moreover, there is
a commutative diagram
$$
\xymatrix{
&&W\ar@{->}[ld]_{\gamma}\ar@{-->}[dr]^{\phi}&\\%
&X\ar@{-->}[rr]_{\psi}&&\mathbb{P}(1,1,3),&}
$$
where $\psi$ is the natural projection and $\phi$ is the rational
map given by $|-3K_{W}|$.

Let $\mathcal{D}$ be the proper transform of $\mathcal{M}$ on $W$.
Then
$$
\mathcal{D}\sim_{\mathbb{Q}}\gamma^{*}(\mathcal{M})-cG-dH,
$$
where $c>r/4$ and $d>r/5$ by \cite{Ka96}.

The natural projection $\psi$ has a one-dimensional family of
fibers that have a singularity at the singular point $O$. Let $C$
be the proper transform on the variety $W$ of a sufficiently
general fiber of the projection $\psi$ that is singular at the
point $O$. Intersecting a general surface of the linear system
$\mathcal{D}$ with the curve $C$, we obtain the inequality
$$
c+2d\leqslant\frac{13r}{20},
$$
which is impossible because $c>r/4$ and $d>r/5$.

Let $S$ be a sufficiently general surface in the linear system
$|-K_{X}|$ and $L$ be the curve on the hypersurface $X$ cut by the
equations $x_1=x_2=0$. Then $S$ is a $K3$ surface whose singular
points are the singular points of the hypersurface $X$. Moreover,
one can easily show that the point $P$ is a singular point of type
$\mathbb{A}_{2}$ on the surface $S$, the point $Q$ is a singular
point of type $\mathbb{A}_{3}$ on the surface $S$, and the point
$O$ is a singular point of type $\mathbb{A}_{4}$ on the surface
$S$. The curve $L$ is a smooth rational curve passing through $P$,
$Q$, and $O$. We have
$$
\mathcal{M}\vert_{S}=\mathcal{P}+\mathrm{mult}_{L}(\mathcal{M})L,
$$
where $\mathcal{P}$ is a linear system on $S$ without fixed
components. Moreover, it immediately follows from \cite{Ka96} that
the inequality $\mathrm{mult}_{L}(\mathcal{M})\leqslant r$ holds
because we already proved that the singularities of $(X,
\frac{1}{r} \mathcal{M})$ are canonical at least at one of the
points $P$, $Q$, and $O$.

Finally, we suppose that the singularities of the log pair $(X,
\frac{1}{r} \mathcal{M})$ are not canonical at the singular points
$Q$ and $P$. Let $\pi:Y\to S$ be the composition of the weighted
blow ups at the points $P$ and $Q$ that are induced by the
Kawamata blow ups of the hypersurface $X$ at the singular points
$P$ and $Q$. Then
$$
\mathcal{H}+\mathrm{mult}_{L}(\mathcal{M})\bar{L}\sim_{\mathbb{Q}}
\pi^{*}\big(-rK_{X}\vert_{S}\big)-m_{1}E_{1}-m_{2}E_{2},
$$
where $E_{1}$ and $E_{2}$ are the $\pi$-exceptional divisors
dominating $P$ and $Q$, respectively, $\bar{L}$ is the proper
transform on the surface $Y$ of the curve $L$, $\mathcal{H}$ is
the proper transform on $Y$ of the linear system $\mathcal{P}$,
$m_{1}$ and $m_{2}$ are rational numbers. Then
$\bar{L}^{2}=-1/30$, but it follows from the paper \cite{Ka96}
that the inequalities $m_{1}>r/3$ and $m_{2}>r/4$ hold.

The curve $\bar{L}$ intersects the curves $E_{1}$ and $E_{2}$ at
singular points of types $\mathbb{A}_{1}$ and $\mathbb{A}_{2}$
respectively. Therefore, the inequalities  $\bar{L}\cdot
E_{1}\geqslant 1/2$ and $\bar{L}\cdot E_{2}\geqslant 1/3$ hold.
Consequently, we obtain
$$
-\frac{r}{30}\leqslant\big(\mathcal{H}+\mathrm{mult}_{L}(\mathcal{M})
\bar{L}\big)\cdot\bar{L}\leqslant \frac{13r}{60}
-\frac{m_{1}}{2}-\frac{m_{2}}{3}<-\frac{r}{30},%
$$
which is a contradiction.
\end{proof}

Therefore, to conclude the proof of
Theorem~\ref{theorem:birational-autmorphisms}, it is enough to
consider the cases
$$
\nn\in\{5, 12, 13, 25, 31, 33, 38, 58\}.
$$
In these  cases  the group $\Gamma_X$ is generated by  two
involutions $\tau_1$ and $\tau_2$. We must show that the group
$\Gamma_X$ is the free product of $\tau_{1}$ and $\tau_{2}$.

Perhaps, the simplest possible way to prove the required claim is
to use the arguments of the proofs of
Lemmas~\ref{lemma:free-group-7} and \ref{lemma:free-group-20}. For
example, the arguments used during the elimination of the points
$Q$ and $O$ in the proof of Lemma~\ref{lemma:free-group-20}
immediately imply the required claim in the case $\nn=5$. However,
we choose an alternative approach.

Let $\psi:X\dasharrow\mathbb{P}(1,a_{2},a_{3})$ be the natural
projection.

\begin{lemma}
\label{lemma:irreducible-fibers-codimension-2}  There are only
finitely many reducible fibers of $\psi$.
\end{lemma}

\begin{proof}
We consider only the case $\nn=58$ because the other cases are
similar. Then $X$ is a sufficiently general hypersurface of degree
$24$ in $\mathbb{P}(1,3,4,7,10)$. It is enough to show that the
fiber $C$ of the projection $\psi$ over a point $(p_1: p_2:
p_3)\in\mathbb{P}(1,3,4)$ is irreducible if $p_1\ne 0$ and
$(p_{1}: p_{2}: p_{3})$ belongs to the complement to a finite set.

By construction, the fiber $C$ is a curve of degree $24/70$ in
$\mathbb{P}(1,7,10)=\mathrm{Proj}(\mathbb{C}[x_1,x_4,x_5])$, where
$\mathrm{wt}(x_{1})=1$, $\mathrm{wt}(x_{4})=7$, and
$\mathrm{wt}(x_{5})=10$. If the curve $C$ is reducible, it must
contain a curve of degree $1/70$, $1/10$, or $1/7$. However, we
have a unique curve of degree $1/70$ in $\mathbb{P}(1,7,10)$,
namely, the curve defined by $x_1=0$. Hence, the fiber $C$ cannot
contain the curve of degree $1/70$ by the generality of the
hypersurface $X$.

Let $\mathcal{X}=|\mathcal{O}_{\mathbb{P}(1,\,3,4,7,10)}(24)|$ and
$\mathcal{C}_{1/7}$ be the set of curves in
$\mathbb{P}(1,3,4,7,10)$ given by
$$
\lambda x_1^3+ x_2=\mu
x_1^4+x_3=\nu_0x_1^{10}+\nu_1x_5+\nu_2x_1^3x_4=0,
$$
where $(\nu_0: \nu_1: \nu_2)\in\mathbb{P}^2$ and $(\lambda,
\mu)\in\mathbb{C}^{2}$. Put
$$
\Gamma=\{ (X, C)\in
\mathcal{X}\times\mathcal{C}_{1/7} \hspace{2mm}|\hspace{2mm} C\subset X\}%
$$
and consider the natural projections $f:\Gamma\to \mathcal{X}$ and
$g:\Gamma\to \mathcal{C}_{1/7}$. Then the pro\-jec\-ti\-on $g$ is
surjective,
$\mathrm{dim}(g^{-1}(x_2=x_3=x_5=0))=\mathrm{dim}(\mathcal{X})-4$,
and $\mathrm{dim}(\mathcal{C}_{1/7})=4$. Thus, we have
$$
\mathrm{dim}(\mathcal{X})\geqslant\mathrm{dim}(\Gamma),
$$
which implies that $X$ contains finitely many curves of degree
$1/7$. Similarly, it is impossible to have infinitely many curves
of degree $1/10$ on $X$. Therefore, the fiber $C$ is irreducible
whenever the point $P$ is in the outside of the finitely many
points in $\mathbb{P}(1,3,4)$ and not in the hyperplane $x_1=0$.
Consequently, the statement for the case $\nn=58$ is true.
\end{proof}

The rational map $\psi$ is not defined at two distinct points of
the hypersurface $X$, which we denote by $P$ and $Q$. Let $C$ be a
very general fiber of the map $\psi$. Then $C$ is a smooth
elliptic curve passing through the points $P$ and $Q$. Moreover,
the following well known result implies that the divisor $P-Q$ is
not a torsion in $\mathrm{Pic}(C)$.

\begin{lemma}
\label{lemma:elliptic-surfaces} Let $\tau:S\to\mathbb{P}^{1}$ be
an elliptic fibration such that the surface $S$ is normal and  all
fibers of the elliptic fibration $\tau$ are irreducible. Suppose
that there are distinct disjoint sections $C_{1}$ and $C_{2}$ of
the elliptic fibration $\tau$ such $C_{1}^{2}<0$ and
$C_{2}^{2}<0$. Then for a very general fiber $L$ of the elliptic
fibration $\tau$ the divisor $(C_{1}-C_{2})\vert_{L}$ is not a
torsion in $\mathrm{Pic}(L)$.
\end{lemma}

\begin{proof}
For every natural number $n$ we have
$$
n\big(C_{1}-C_{2}\big)\vert_{L}\sim 0\Rightarrow C_{1}-C_{2}\equiv\Sigma, %
$$
where $\Sigma$ is a $\mathbb{Q}$-divisor on the surface $S$ whose
support is contained in the fibers of the elliptic fibration
$\tau$. On the other hand, because all fibers of $\tau$ are
irreducible, the curves $C_{1}$, $C_{2}$, and $L$ are linearly
dependent in the group
$\mathrm{Div}(S)\otimes\mathbb{Q}\slash\equiv$. However,
$$
\left|\begin{array}{ccc}%
C_{1}^{2}        &C_{1}\cdot C_{2} & C_{1}\cdot L \\ %
C_{1}\cdot C_{2} &C_{2}^{2}        & C_{2}\cdot L \\ %
C_{1}\cdot L     &C_{2}\cdot L     & L^{2}        \\ %
\end{array}\right|=-C_{1}^{2}-C_{2}^{2}\ne 0,
$$
which contradicts  the linear dependence of the curves $C_{1}$,
$C_{2}$, and $L$.
\end{proof}

The curve $C$ is invariant under the action of the birational
involutions $\tau_{1}$ and $\tau_{2}$. Moreover, up to relabelling
the involutions $\tau_{1}$ and $\tau_{2}$ act on the elliptic
curve $C$ by reflections with respect to the points $Q$ and $P$,
respectively. Hence, the composition $\tau_{1}\circ\tau_{2}$ acts
the smooth elliptic curve $C$ by the translation on $2(P-Q)$.
Therefore, the composition $(\tau_{1}\circ\tau_{2})^{n}$ never
acts identically on the curve $C$ for any natural number $n\ne 0$
because the divisor $P-Q$ on the curve $C$ is not a torsion in
$\mathrm{Pic}(C)$. Hence, the group $\Gamma_X$ is the free
products of the birational involutions $\tau_{1}$ and $\tau_{2}$,
which concludes the proof of
Theorem~\ref{theorem:birational-autmorphisms}.

\section{Potential density.}
\label{section:simple-cases}

 Suppose that the hypersurface
$X$ is defined over a number field $\FF$. The purpose of this
section is to  complete the proof of
Proposition~\ref{proposition:potential-density} by proving the
potential density of the set of rational points of the
hypersurface $X$ in the cases $\nn=11$ and $19$.

\begin{lemma}
\label{lemma:potential-density-19} Suppose that $\nn=19$. Then
rational points on $X$ are potentially dense.
\end{lemma}

\begin{proof}
The 3-fold $X$ is a general hypersurface in
$\mathbb{P}_\FF(1,2,3,3,4)$ given by the equation
$$
\sum_{{i, j, k, l, m\geq
0}\atop{i+a_{2}j+a_{3}k+a_{4}l+a_{5}m=12}}a_{ijklm}
x_1^{i}x_2^{j}x_3^{k}x_4^{l}x_5^{m}=0,%
$$
where $a_{ijklm}\in\mathbb{F}$ and we may assume that
$a_{00040}=0$ and $a_{00003}=1$ possibly after replacing the field
$\mathbb{F}$ by its finite extension. Let $P=(0:0:0:1:0)$. Then
$X$ has a cyclic quotient singularity of type
${\frac{1}{3}}(1,2,1)$ at the point $P$.

Let $\alpha:V\to X$ be the Kawamata blow up at $P$. Then the
equality $-K_{V}^{3}=0$ holds, the linear system $|-6K_{V}|$ has
no base point, and
$$
K_{V}\sim_{\mathbb{Q}}\alpha^{*}(K_{X})+{\frac{1}{3}}E,%
$$
where $E=\alpha^{-1}(P)\cong\mathbb{P}(1,1,2)$. Let
$\psi:V\to\mathbb{P}(1,2,3)$ be the morphism given by the linear
system $|-6K_{V}|$. Then $\psi$ is an elliptic fibration (see the
proof of
Lemma~\ref{lemma:existence-of-elliptic-fibration-7-11-19}).

The restriction $\psi\vert_{E}:E\to\mathbb{P}(1,2,3)$ is a triple
cover,  namely, the divisor $E$ is a $3$-section of the elliptic
fibration $\psi$. In the case when $\psi\vert_{E}$ is branched at
a point contained in a smooth fiber of $\psi$, the set of rational
points on $V$ is potentially dense (see \cite{BoTsch99}) because
$E$ is a rational surface. Therefore, it is enough to find a
smooth fiber $C$ of the fibration $\psi$ such that the
intersection $C\cap E$ consists of at most two points.

Let $Z$ be the curve on $X$ given by the equations $x_2=\lambda
x_1^{2}$ and $x_3=\mu x_1^{3}$, where $\lambda,\mu\in\mathbb{F}$,
and $\hat{Z}=\alpha^{-1}(Z)$. Then $\hat{Z}$ is a fiber of $\psi$.
The intersection $\hat{Z}\cap E$ consists of three different
points if and only if $Z$ has an ordinary triple point at $P$.
However, the curve $Z$ has an ordinary triple point at the point
$P$ if and only if the homogeneous polynomial
$$
f(x_1,x_5)=x_5^{3}+a_{10012}x_1x_5^{2}+x_5x_1^{2}(a_{20021}+
\lambda a_{01021})+x_1^{3}(\mu a_{00130}+\lambda a_{10030}+a_{30030})%
$$
has three distinct roots. Now if we  put
$$
\lambda=\frac{a_{10012}^2-4a_{20021}}{4a_{01021}}\ \text{and}\ \mu=-{\frac{\lambda a_{10030}+a_{30030}}{a_{00130}}},%
$$
then the generality of the hypersurface $X$ together with the
Bertini theorem implies that the curve $\hat{Z}$ is smooth but the
intersection $\hat{Z}\cap E$ consists of only two different
points.
\end{proof}

To prove the potential density of the case $\nn=11$, we first
consider a general surface in $|-K_X|$.

\begin{lemma}
\label{lemma:potential-density-K3-11} Let $Y$ be a general surface
in $|-K_{X}|$. Suppose that at least one singular point of $Y$ is
defined over the field $\mathbb{F}$. Then the set of
$\mathbb{F}$-rational points of $Y$ is Zariski dense.
\end{lemma}

\begin{proof}
We have a hypersurface $Y\subset\mathbb{P}(1,2,2,5)$ which can be
given by the equation
$$
x_4^{2}=x_1^{2}f_{4}(x_2,x_3)+x_1^{4}f_{3}(x_2,x_3)+x_1^{6}f_{2}(x_2,x_3)+x_1^{8}f_{1}(x_2,x_3)+x_1^{10}+x_3g_{4}(x_2,x_3),
$$
where $f_{i}$ and $g_{i}$ are general homogeneous polynomials of
degree $i$.

Let $P$ be the point $(0:1:0:0)$ and $\mathcal{H}$ be the pencil
of curves on $Y$ given by the equations $\lambda x_1^{2}+\mu
x_3=0$, $(\lambda:\mu)\in\mathbb{P}_\FF^{1}$. Then $Y$ has a
singularity of type $\mathbb{A}_{1}$ at the point $P$ which is a
unique base point of the pencil $\mathcal{H}$.

Let $C$ be the curve in $\mathcal{H}$ corresponding to the point
$(\lambda: \mu)\in\mathbb{P}_\FF^{1}$ and
$$
f_{4}(x_2,x_3)=\sum_{i=0}^{4}\alpha_{i}x_2^{i}x_3^{4-i},\ g_{4}(x_2,x_3)=\sum_{i=0}^{4}\beta_{i}x_2^{i}x_3^{4-i},%
$$
where $\alpha_{i}$ and $\beta_{i}$ are sufficiently general
constants. Then the curve $C$ has an ordinary double point at the
point $P$ when $(\lambda: \mu)\ne(1:0)$ and $(\lambda:
\mu)\ne(\alpha_{4}:\beta_{4})$. Let $F$ be the curve in the pencil
$\mathcal{H}$ corresponding to the point $(\lambda:
\mu)=(\alpha_{4}:\beta_{4})$ and $L$ be the curve on the surface
$Y$ given by the equation $x_1=0$. Then $F$ is smooth in the
outside of $P$ and has an ordinary cusp at $P$, while $L$ is a
smooth rational curve.

Let $\pi:W\to Y$ be the blow up at the point  $P$, $E$ be the
$\pi$-exceptional divisor, and $\mathcal{B}$ be the proper
transform of the pencil $\mathcal{H}$ on the surface $W$. Then
$\mathcal{B}$ has no base point and induces an elliptic fibration
$\psi:W\to\mathbb{P}^{1}$. The proper transform $\hat{F}$ of $F$
by $\pi$ is a smooth elliptic fiber of the fibration $\psi$.
Moreover, the restriction $\pi\vert_{E}:E\to\mathbb{P}^{1}$ is a
double cover branched at the point  $\hat{F}\cap E$. Because  the
set of all $\mathbb{F}$-rational points of the curve $E$ is
Zariski dense, it follows from \cite{BoTsch99} that the set of
$\mathbb{F}$-rational points of the surface $S$ is Zariski dense.
\end{proof}

Because we may assume that the singular points of $X$ are
$\FF$-rational by replacing $\FF$ by its finite extension, one can
easily prove the density of $\FF$-rational points on $X$ with the
lemma above.

\section{Appendix}
\label{section:weighted-Fanos}

The list of quasismooth anticanonically embedded weighted Fano
3-fold hypersurfaces is found in \cite{IF00}. The completeness of
the list is proved in \cite{JoKo01}.

\Small{
\renewcommand\arraystretch{1.4}
\singlespacing
\begin{longtable}{|c|c|c|c|c|c|c|c|c|}
\caption{Weighted Fano hypersurfaces of degree $d$ in
$\mathbb{P}(1,a_{2},a_{3},a_{4},a_{5})$.\label{table:weighted-Fanos}
}
\\ \hline
$\nn$ & $d$ & $a_{2}$ & $a_{3}$ & $a_{4}$ & $a_{5}$ & $-K_{X}^{3}$ & $\mathrm{Sing}(X)$ &  $\Gamma_X$ \\

\endfirsthead
\hline \multicolumn{9}{|l|}{\scriptsize\slshape continued from
previous page}\\ \hline $\nn$  & $d$ & $a_{2}$ & $a_{3}$ & $a_{4}$
& $a_{5}$ & $-K_{X}^{3}$ & $\mathrm{Sing}(X)$ &  $\Gamma_X$ \\
\hline \hline\endhead \hline
\multicolumn{9}{|l|}{\scriptsize\slshape continued on next page}\\
\hline
\endfoot
\hline

\endlastfoot

\hline\hline
$1$  & $4$ & $1$ & $1$ & $1$ & $1$ & $4$ & $\varnothing$& \FGZero \\%
\hline
$2$  & $5$ & $1$ & $1$ & $1$ & $2$ & $5/2$ & ${\frac{1}{2}(1,1,1)}$& \FGOne\\%
\hline
$3$  & $6$ & $1$ & $1$ & $1$ & $3$ & $2$ & $\varnothing$& \FGZero\\%
\hline
$4$  & $6$ & $1$ & $1$ & $2$ & $2$ & $3/2$ & ${3\times\frac{1}{2}(1,1,2)}$& \EG\\%
\hline
$5$  & $7$ & $1$ & $1$ & $2$ & $3$ & $7/6$ & ${\frac{1}{2}(1,1,1),\frac{1}{3}(1,1,2)}$& \FGTwo\\%
\hline
$6$  & $8$ & $1$ & $1$ & $2$ & $4$ & $1$ & ${2\times\frac{1}{2}(1,1,1)}$& \FGTwo \\%
\hline
$7$  & $8$ & $1$ & $2$ & $2$ & $3$ & $2/3$ & ${4\times\frac{1}{2}(1,1,1),\frac{1}{3}(1,1,2)}$& \FGFive\\%
\hline
$8$  & $9$ & $1$ & $1$ & $3$ & $4$ & $3/4$ & ${\frac{1}{4}(1,1,3)}$& \FGOne\\%
\hline
$9$  & $9$ & $1$ & $2$ & $3$ & $3$ & $1/2$ & ${\frac{1}{2}(1,1,1), 3\times\frac{1}{3}(1,1,2)}$& \EG\\%
\hline
$10$  & $10$ & $1$ & $1$ & $3$ & $5$ & $2/3$ & ${\frac{1}{3}(1,1,2)}$& \FGZero\\%
\hline
$11$  & $10$ & $1$ & $2$ & $2$ & $5$ & $1/2$ & ${5\times\frac{1}{2}(1,1,1)}$& \FGZero\\%
\hline
$12$  & $10$ & $1$ & $2$ & $3$ & $4$ & $5/12$ & ${2\times\frac{1}{2}(1,1,1),\frac{1}{3}(1,1,2),\frac{1}{4}(1,1,3)}$& \FGTwo\\%
\hline
$13$  & $11$ & $1$ & $2$ & $3$ & $5$ & $11/30$ & ${\frac{1}{2}(1,1,1),\frac{1}{3}(1,1,2),\frac{1}{5}(1,2,3)}$& \FGTwo\\%
\hline
$14$  & $12$ & $1$ & $1$ & $4$ & $6$ & $1/2$ & ${\frac{1}{2}(1,1,1)}$& \FGZero\\%
\hline
$15$  & $12$ & $1$ & $2$ & $3$ & $6$ & $1/3$ & ${2\times\frac{1}{2}(1,1,1), 2\times\frac{1}{3}(1,1,2)}$& \FGTwo\\%
\hline
$16$  & $12$ & $1$ & $2$ & $4$ & $5$ & $3/10$ & ${3\times\frac{1}{2}(1,1,1),\frac{1}{5}(1,1,4)}$& \FGOne\\%
\hline
$17$  & $12$ & $1$ & $3$ & $4$ & $4$ & $1/4$ & ${3\times\frac{1}{4}(1,1,3)}$& \EG\\%
\hline
$18$  & $12$ & $2$ & $2$ & $3$ & $5$ & $1/5$ & ${6\times\frac{1}{2}(1,1,1),\frac{1}{5}(1,2,3)}$& \FGOne\\%
\hline
$19$  & $12$ & $2$ & $3$ & $3$ & $4$ & $1/6$ & ${3\times\frac{1}{2}(1,1,1), 4\times\frac{1}{3}(1,1,2)}$& \FGZero\\%
\hline
$20$  & $13$ & $1$ & $3$ & $4$ & $5$ & $13/60$ & ${\frac{1}{3}(1,1,2),\frac{1}{4}(1,1,3),\frac{1}{5}(1,1,4)}$& \FGThree\\%
\hline
$21$  & $14$ & $1$ & $2$ & $4$ & $7$ & $1/4$ & ${3\times\frac{1}{2}(1,1,1),\frac{1}{4}(1,1,3)}$& \FGZero\\%
\hline
$22$  & $14$ & $2$ & $2$ & $3$ & $7$ & $1/6$ & ${7\times\frac{1}{2}(1,1,1),\frac{1}{3}(1,1,2)}$& \FGZero\\%
\hline
$23$  & $14$ & $2$ & $3$ & $4$ & $5$ & $7/60$ & ${3\times\frac{1}{2}(1,1,1),\frac{1}{3}(1,1,2),\frac{1}{4}(1,1,3),\frac{1}{5}(1,2,3)}$& \FGTwo\\%
\hline
$24$  & $15$ & $1$ & $2$ & $5$ & $7$ & $3/14$ & ${\frac{1}{2}(1,1,1),\frac{1}{7}(1,2,5)}$& \FGOne\\%
\hline
$25$  & $15$ & $1$ & $3$ & $4$ & $7$ & $5/28$ & ${\frac{1}{4}(1,1,3),\frac{1}{7}(1,3,4)}$& \FGTwo\\%
\hline
$26$  & $15$ & $1$ & $3$ & $5$ & $6$ & $1/6$ & ${2\times\frac{1}{3}(1,1,2),\frac{1}{6}(1,1,5)}$& \FGOne\\%
\hline
$27$  & $15$ & $2$ & $3$ & $5$ & $5$ & $1/10$ & ${\frac{1}{2}(1,1,1), 3\times\frac{1}{5}(1,2,3)}$& \EG\\%
\hline
$28$  & $15$ & $3$ & $3$ & $4$ & $5$ & $1/12$ & ${5\times\frac{1}{3}(1,1,2),\frac{1}{4}(1,1,3)}$& \FGZero\\%
\hline
$29$  & $16$ & $1$ & $2$ & $5$ & $8$ & $1/5$ & ${2\times\frac{1}{2}(1,1,1),\frac{1}{5}(1,2,3)}$& \FGZero\\%
\hline
$30$  & $16$ & $1$ & $3$ & $4$ & $8$ & $1/6$ & ${\frac{1}{3}(1,1,2), 2\times\frac{1}{4}(1,1,3)}$& \FGTwo\\%
\hline
$31$  & $16$ & $1$ & $4$ & $5$ & $6$ & $2/15$ & ${\frac{1}{2}(1,1,1),\frac{1}{5}(1,1,4),\frac{1}{6}(1,1,5)}$& \FGTwo\\%
\hline
$32$  & $16$ & $2$ & $3$ & $4$ & $7$ & $2/21$ & ${4\times\frac{1}{2}(1,1,1),\frac{1}{3}(1,1,2),\frac{1}{7}(1,3,4)}$& \FGOne\\%
\hline
$33$  & $17$ & $2$ & $3$ & $5$ & $7$ & $17/210$ & ${\frac{1}{2}(1,1,1),\frac{1}{3}(1,1,2),\frac{1}{5}(1,2,3),\frac{1}{7}(1,2,5)}$& \FGTwo\\%
\hline
$34$  & $18$ & $1$ & $2$ & $6$ & $9$ & $1/6$ & ${3\times\frac{1}{2}(1,1,1),\frac{1}{3}(1,1,2)}$& \FGZero\\%
\hline
$35$  & $18$ & $1$ & $3$ & $5$ & $9$ & $2/15$ & ${2\times\frac{1}{3}(1,1,2),\frac{1}{5}(1,1,4)}$& \FGZero\\%
\hline
$36$  & $18$ & $1$ & $4$ & $6$ & $7$ & $3/28$ & ${\frac{1}{4}(1,1,3),\frac{1}{2}(1,1,1),\frac{1}{7}(1,1,6)}$& \FGTwo\\%
\hline
$37$  & $18$ & $2$ & $3$ & $4$ & $9$ & $1/12$ & ${4\times\frac{1}{2}(1,1,1), 2\times\frac{1}{3}(1,1,2),\frac{1}{4}(1,1,3)}$& \FGZero\\%
\hline
$38$  & $18$ & $2$ & $3$ & $5$ & $8$ & $3/40$ & ${2\times\frac{1}{2}(1,1,1),\frac{1}{5}(1,2,3),\frac{1}{8}(1,3,5)}$& \FGTwo\\%
\hline
$39$  & $18$ & $3$ & $4$ & $5$ & $6$ & $1/20$ & ${3\times\frac{1}{3}(1,1,2),\frac{1}{4}(1,1,3),\frac{1}{2}(1,1,1),\frac{1}{5}(1,1,4)}$& \FGZero\\%
\hline
$40$  & $19$ & $3$ & $4$ & $5$ & $7$ & $19/420$ & ${\frac{1}{3}(1,1,2),\frac{1}{4}(1,1,3),\frac{1}{5}(1,2,3),\frac{1}{7}(1,3,4)}$& \FGTwo\\%
\hline
$41$  & $20$ & $1$ & $4$ & $5$ & $10$ & $1/10$ & ${\frac{1}{2}(1,1,1), 2\times\frac{1}{5}(1,1,4)}$& \FGTwo\\%
\hline
$42$  & $20$ & $2$ & $3$ & $5$ & $10$ & $1/15$ & ${2\times\frac{1}{2}(1,1,1),\frac{1}{3}(1,1,2), 2\times\frac{1}{5}(1,2,3)}$& \FGTwo\\%
\hline
$43$  & $20$ & $2$ & $4$ & $5$ & $9$ & $1/18$ & ${5\times\frac{1}{2}(1,1,1),\frac{1}{9}(1,4,5)}$& \FGOne\\%
\hline
$44$  & $20$ & $2$ & $5$ & $6$ & $7$ & $1/21$ & ${3\times\frac{1}{2}(1,1,1),\frac{1}{6}(1,1,5),\frac{1}{7}(1,2,5)}$& \FGTwo\\%
\hline
$45$  & $20$ & $3$ & $4$ & $5$ & $8$ & $1/24$ & ${\frac{1}{3}(1,1,2), 2\times\frac{1}{4}(1,1,3),\frac{1}{8}(1,3,5)}$& \FGOne\\%
\hline
$46$  & $21$ & $1$ & $3$ & $7$ & $10$ & $1/10$ & ${\frac{1}{10}(1,3,7)}$& \FGOne\\%
\hline
$47$  & $21$ & $1$ & $5$ & $7$ & $8$ & $3/40$ & ${\frac{1}{5}(1,2,3),\frac{1}{8}(1,1,7)}$& \FGOne\\%
\hline
$48$  & $21$ & $2$ & $3$ & $7$ & $9$ & $1/18$ & ${\frac{1}{2}(1,1,1), 2\times\frac{1}{3}(1,1,2),\frac{1}{9}(1,2,7)}$& \FGOne\\%
\hline
$49$  & $21$ & $3$ & $5$ & $6$ & $7$ & $1/30$ & ${3\times\frac{1}{3}(1,1,2),\frac{1}{5}(1,2,3),\frac{1}{6}(1,1,5)}$& \FGZero\\%
\hline
$50$  & $22$ & $1$ & $3$ & $7$ & $11$ & $2/21$ & ${\frac{1}{3}(1,1,2),\frac{1}{7}(1,3,4)}$& \FGZero\\%
\hline
$51$  & $22$ & $1$ & $4$ & $6$ & $11$ & $1/12$ & ${\frac{1}{4}(1,1,3),\frac{1}{2}(1,1,1),\frac{1}{6}(1,1,5)}$& \FGZero\\%
\hline
$52$  & $22$ & $2$ & $4$ & $5$ & $11$ & $1/20$ & ${5\times\frac{1}{2}(1,1,1),\frac{1}{4}(1,1,3),\frac{1}{5}(1,1,4)}$& \FGZero\\%
\hline
$53$  & $24$ & $1$ & $3$ & $8$ & $12$ & $1/12$ & ${2\times\frac{1}{3}(1,1,2),\frac{1}{4}(1,1,3)}$& \FGZero\\%
\hline
$54$  & $24$ & $1$ & $6$ & $8$ & $9$ & $1/18$ & ${\frac{1}{2}(1,1,1),\frac{1}{3}(1,1,2),\frac{1}{9}(1,1,8)}$& \FGOne\\%
\hline
$55$  & $24$ & $2$ & $3$ & $7$ & $12$ & $1/21$ & ${2\times\frac{1}{2}(1,1,1), 2\times\frac{1}{3}(1,1,2),\frac{1}{7}(1,2,5)}$& \FGZero\\%
\hline
$56$  & $24$ & $2$ & $3$ & $8$ & $11$ & $1/22$ & ${3\times\frac{1}{2}(1,1,1),\frac{1}{11}(1,3,8)}$& \FGOne\\%
\hline
$57$  & $24$ & $3$ & $4$ & $5$ & $12$ & $1/30$ & ${2\times\frac{1}{3}(1,1,2), 2\times\frac{1}{4}(1,1,3),\frac{1}{5}(1,2,3)}$& \FGZero\\%
\hline
$58$  & $24$ & $3$ & $4$ & $7$ & $10$ & $1/35$ & ${\frac{1}{2}(1,1,1),\frac{1}{7}(1,3,4),\frac{1}{10}(1,3,7)}$& \FGTwo\\%
\hline
$59$  & $24$ & $3$ & $6$ & $7$ & $8$ & $1/42$ & ${4\times\frac{1}{3}(1,1,2),\frac{1}{2}(1,1,1),\frac{1}{7}(1,1,6)}$& \FGZero\\%
\hline
$60$  & $24$ & $4$ & $5$ & $6$ & $9$ & $1/45$ & ${2\times\frac{1}{2}(1,1,1),\frac{1}{5}(1,1,4),\frac{1}{3}(1,1,2),\frac{1}{9}(1,4,5)}$& \FGOne\\%
\hline
$61$  & $25$ & $4$ & $5$ & $7$ & $9$ & $5/252$ & ${\frac{1}{4}(1,1,3),\frac{1}{7}(1,2,5),\frac{1}{9}(1,4,5)}$& \FGTwo\\%
\hline
$62$  & $26$ & $1$ & $5$ & $7$ & $13$ & $2/35$ & ${\frac{1}{5}(1,2,3),\frac{1}{7}(1,1,6)}$& \FGZero\\%
\hline
$63$  & $26$ & $2$ & $3$ & $8$ & $13$ & $1/24$ & ${3\times\frac{1}{2}(1,1,1),\frac{1}{3}(1,1,2),\frac{1}{8}(1,3,5)}$& \FGZero\\%
\hline
$64$  & $26$ & $2$ & $5$ & $6$ & $13$ & $1/30$ & ${4\times\frac{1}{2}(1,1,1),\frac{1}{5}(1,2,3),\frac{1}{6}(1,1,5)}$& \FGZero\\%
\hline
$65$  & $27$ & $2$ & $5$ & $9$ & $11$ & $3/110$ & ${\frac{1}{2}(1,1,1),\frac{1}{5}(1,1,4),\frac{1}{11}(1,2,9)}$& \FGOne\\%
\hline
$66$  & $27$ & $5$ & $6$ & $7$ & $9$ & $1/70$ & ${\frac{1}{5}(1,1,4),\frac{1}{6}(1,1,5),\frac{1}{3}(1,1,2),\frac{1}{7}(1,2,5)}$& \FGZero\\%
\hline
$67$  & $28$ & $1$ & $4$ & $9$ & $14$ & $1/18$ & ${\frac{1}{2}(1,1,1),\frac{1}{9}(1,4,5)}$& \FGZero\\%
\hline
$68$  & $28$ & $3$ & $4$ & $7$ & $14$ & $1/42$ & ${\frac{1}{3}(1,1,2),\frac{1}{2}(1,1,1), 2\times\frac{1}{7}(1,3,4)}$& \FGTwo\\%
\hline
$69$  & $28$ & $4$ & $6$ & $7$ & $11$ & $1/66$ & ${2\times\frac{1}{2}(1,1,1),\frac{1}{6}(1,1,5),\frac{1}{11}(1,4,7)}$& \FGOne\\%
\hline
$70$  & $30$ & $1$ & $4$ & $10$ & $15$ & $1/20$ & ${\frac{1}{4}(1,1,3),\frac{1}{2}(1,1,1),\frac{1}{5}(1,1,4)}$& \FGZero\\%
\hline
$71$  & $30$ & $1$ & $6$ & $8$ & $15$ & $1/24$ & ${\frac{1}{2}(1,1,1),\frac{1}{3}(1,1,2),\frac{1}{8}(1,1,7)}$& \FGZero\\%
\hline
$72$  & $30$ & $2$ & $3$ & $10$ & $15$ & $1/30$ & ${3\times\frac{1}{2}(1,1,1), 2\times\frac{1}{3}(1,1,2),\frac{1}{5}(1,2,3)}$& \FGZero\\%
\hline
$73$  & $30$ & $2$ & $6$ & $7$ & $15$ & $1/42$ & ${5\times\frac{1}{2}(1,1,1),\frac{1}{3}(1,1,2),\frac{1}{7}(1,1,6)}$& \FGZero\\%
\hline
$74$  & $30$ & $3$ & $4$ & $10$ & $13$ & $1/52$ & ${\frac{1}{4}(1,1,3),\frac{1}{2}(1,1,1),\frac{1}{13}(1,3,10)}$& \FGOne\\%
\hline
$75$  & $30$ & $4$ & $5$ & $6$ & $15$ & $1/60$ & ${\frac{1}{4}(1,1,3), 2\times\frac{1}{2}(1,1,1), 2\times\frac{1}{5}(1,1,4),\frac{1}{3}(1,1,2)}$& \FGZero\\%
\hline
$76$  & $30$ & $5$ & $6$ & $8$ & $11$ & $1/88$ & ${\frac{1}{2}(1,1,1),\frac{1}{8}(1,3,5),\frac{1}{11}(1,5,6)}$& \FGTwo\\%
\hline
$77$  & $32$ & $2$ & $5$ & $9$ & $16$ & $1/45$ & ${2\times\frac{1}{2}(1,1,1),\frac{1}{5}(1,1,4),\frac{1}{9}(1,2,7)}$& \FGZero\\%
\hline
$78$  & $32$ & $4$ & $5$ & $7$ & $16$ & $1/70$ & ${2\times\frac{1}{4}(1,1,3),\frac{1}{5}(1,1,4),\frac{1}{7}(1,5,2)}$& \FGZero\\%
\hline
$79$  & $33$ & $3$ & $5$ & $11$ & $14$ & $1/70$ & ${\frac{1}{5}(1,1,4),\frac{1}{14}(1,3,11)}$& \FGOne\\%
\hline
$80$  & $34$ & $3$ & $4$ & $10$ & $17$ & $1/60$ & ${\frac{1}{3}(1,1,2),\frac{1}{4}(1,1,3),\frac{1}{2}(1,1,1),\frac{1}{10}(1,3,7)}$& \FGZero\\%
\hline
$81$  & $34$ & $4$ & $6$ & $7$ & $17$ & $1/84$ & ${\frac{1}{4}(1,1,3), 2\times\frac{1}{2}(1,1,1),\frac{1}{6}(1,1,5),\frac{1}{7}(1,4,3)}$& \FGZero\\%
\hline
$82$  & $36$ & $1$ & $5$ & $12$ & $18$ & $1/30$ & ${\frac{1}{5}(1,2,3),\frac{1}{6}(1,1,5)}$& \FGZero\\%
\hline
$83$  & $36$ & $3$ & $4$ & $11$ & $18$ & $1/66$ & ${2\times\frac{1}{3}(1,1,2),\frac{1}{2}(1,1,1),\frac{1}{11}(1,4,7)}$& \FGZero\\%
\hline
$84$  & $36$ & $7$ & $8$ & $9$ & $12$ & $1/168$ & ${\frac{1}{7}(1,2,5),\frac{1}{8}(1,1,7),\frac{1}{4}(1,1,3),\frac{1}{3}(1,1,2)}$& \FGZero\\%
\hline
$85$  & $38$ & $3$ & $5$ & $11$ & $19$ & $2/165$ & ${\frac{1}{3}(1,1,2),\frac{1}{5}(1,1,4),\frac{1}{11}(1,3,8)}$& \FGZero\\%
\hline
$86$  & $38$ & $5$ & $6$ & $8$ & $19$ & $1/120$ & ${\frac{1}{5}(1,1,4),\frac{1}{6}(1,1,5),\frac{1}{2}(1,1,1),\frac{1}{8}(1,3,5)}$& \FGZero\\%
\hline
$87$  & $40$ & $5$ & $7$ & $8$ & $20$ & $1/140$ & ${2\times\frac{1}{5}(1,2,3),\frac{1}{7}(1,1,6),\frac{1}{4}(1,1,3)}$& \FGZero\\%
\hline
$88$  & $42$ & $1$ & $6$ & $14$ & $21$ & $1/42$ & ${\frac{1}{2}(1,1,1),\frac{1}{3}(1,1,2),\frac{1}{7}(1,1,6)}$& \FGZero\\%
\hline
$89$  & $42$ & $2$ & $5$ & $14$ & $21$ & $1/70$ & ${3\times\frac{1}{2}(1,1,1),\frac{1}{5}(1,1,4),\frac{1}{7}(1,2,5)}$& \FGZero\\%
\hline
$90$  & $42$ & $3$ & $4$ & $14$ & $21$ & $1/84$ & ${2\times\frac{1}{3}(1,1,2),\frac{1}{4}(1,1,3),\frac{1}{2}(1,1,1),\frac{1}{7}(1,3,4)}$& \FGZero\\%
\hline
$91$  & $44$ & $4$ & $5$ & $13$ & $22$ & $1/130$ & ${\frac{1}{2}(1,1,1),\frac{1}{5}(1,2,3),\frac{1}{13}(1,4,9)}$& \FGZero\\%
\hline
$92$  & $48$ & $3$ & $5$ & $16$ & $24$ & $1/120$ & ${2\times\frac{1}{3}(1,1,2),\frac{1}{5}(1,1,4),\frac{1}{8}(1,3,5)}$& \FGZero\\%
\hline
$93$  & $50$ & $7$ & $8$ & $10$ & $25$ & $1/280$ & ${\frac{1}{7}(1,3,4),\frac{1}{8}(1,1,7),\frac{1}{2}(1,1,1),\frac{1}{5}(1,2,3)}$& \FGZero\\%
\hline
$94$  & $54$ & $4$ & $5$ & $18$ & $27$ & $1/180$ & ${\frac{1}{4}(1,1,3),\frac{1}{2}(1,1,1),\frac{1}{5}(1,2,3),\frac{1}{9}(1,4,5)}$& \FGZero\\%
\hline
$95$  & $66$ & $5$ & $6$ & $22$ & $33$ & $1/330$ & ${\frac{1}{5}(1,2,3),\frac{1}{2}(1,1,1),\frac{1}{3}(1,1,2),\frac{1}{11}(1,5,6)}$& \FGZero\\%
\hline
\end{longtable}} 
\end{document}